\documentclass[times,preprint]{elsarticle}

\usepackage{lineno}
 \RequirePackage{geometry}
 \global\bibsep=0pt

\usepackage{framed,multirow}

\usepackage{latexsym}

\usepackage{url}
\usepackage{xcolor}
\definecolor{newcolor}{rgb}{.8,.349,.1}

\journal{}

\usepackage{graphicx}
\usepackage{lipsum}
\usepackage{amsfonts}
\usepackage{graphicx}
\usepackage{epstopdf}
\usepackage{algorithmic}
\ifpdf
\DeclareGraphicsExtensions{.eps,.pdf,.png,.jpg}
\else
\DeclareGraphicsExtensions{.eps}
\fi

\usepackage{multicol}        
\usepackage[bottom]{footmisc}

\usepackage{amssymb}
\usepackage{amsmath}
\usepackage{amsthm}
\usepackage{wasysym}

\usepackage{bm}
\usepackage{xspace}

\usepackage{booktabs}
\usepackage{tabularx}

\usepackage{caption}
\usepackage{subfig}

\usepackage{color}
\usepackage[normalem]{ulem}

\usepackage{hyperref}
\usepackage[capitalize]{cleveref}

\usepackage{xargs}  

\renewcommand{\vec}[1]{{\bf #1}}

\def\@thmcountersep{.}

\newcounter{theorem}

\newtheorem{mytheorem}[theorem]{Theorem}

\newtheorem{myremark}[theorem]{Remark}

\newcolumntype{L}[1]{>{\raggedright\arraybackslash}p{#1}} 
\newcolumntype{C}[1]{>{\centering\arraybackslash}p{#1}} 
\newcolumntype{R}[1]{>{\raggedleft\arraybackslash}p{#1}} 
\newcolumntype{Y}{>{\centering\arraybackslash}X} 
\newcolumntype{Z}{>{\raggedleft\arraybackslash}X} 

\numberwithin{theorem}{section}

\def\@thmcountersep{.}

\makeatletter
\newcommand{\manuallabel}[2]{\def\@currentlabel{#2}\label{#1}}
\makeatother

\newcounter{subcount}

\newcommand{\reva}[1]{#1\xspace}
\newcommand{\revb}[1]{#1\xspace}

\DeclareMathAlphabet{\mathitbf}{OML}{cmm}{b}{it}

\makeatletter
\newcommand\RedeclareMathOperator{%
	\@ifstar{\def\rmo@s{m}\rmo@redeclare}{\def\rmo@s{o}\rmo@redeclare}%
}
\newcommand\rmo@redeclare[2]{%
	\begingroup \escapechar\m@ne\xdef\@gtempa{{\string#1}}\endgroup
	\expandafter\@ifundefined\@gtempa
	{\@latex@error{\noexpand#1undefined}\@ehc}%
	\relax
	\expandafter\rmo@declmathop\rmo@s{#1}{#2}}
\newcommand\rmo@declmathop[3]{%
	\DeclareRobustCommand{#2}{\qopname\newmcodes@#1{#3}}%
}
\@onlypreamble\RedeclareMathOperator
\makeatother

\DeclareMathSymbol{\alpha}{\mathalpha}{letters}{"0B}
\DeclareMathSymbol{\beta}{\mathalpha}{letters}{"0C}
\DeclareMathSymbol{\gamma}{\mathalpha}{letters}{"0D}
\DeclareMathSymbol{\delta}{\mathalpha}{letters}{"0E}
\DeclareMathSymbol{\epsilon}{\mathalpha}{letters}{"0F}
\DeclareMathSymbol{\zeta}{\mathalpha}{letters}{"10}
\DeclareMathSymbol{\eta}{\mathalpha}{letters}{"11}
\DeclareMathSymbol{\theta}{\mathalpha}{letters}{"12}
\DeclareMathSymbol{\iota}{\mathalpha}{letters}{"13}
\DeclareMathSymbol{\kappa}{\mathalpha}{letters}{"14}
\DeclareMathSymbol{\lambda}{\mathalpha}{letters}{"15}
\DeclareMathSymbol{\mu}{\mathalpha}{letters}{"16}
\DeclareMathSymbol{\nu}{\mathalpha}{letters}{"17}
\DeclareMathSymbol{\xi}{\mathalpha}{letters}{"18}
\DeclareMathSymbol{\pi}{\mathalpha}{letters}{"19}
\DeclareMathSymbol{\rho}{\mathalpha}{letters}{"1A}
\DeclareMathSymbol{\sigma}{\mathalpha}{letters}{"1B}
\DeclareMathSymbol{\tau}{\mathalpha}{letters}{"1C}
\DeclareMathSymbol{\upsilon}{\mathalpha}{letters}{"1D}
\DeclareMathSymbol{\phi}{\mathalpha}{letters}{"1E}
\DeclareMathSymbol{\chi}{\mathalpha}{letters}{"1F}
\DeclareMathSymbol{\psi}{\mathalpha}{letters}{"20}
\DeclareMathSymbol{\omega}{\mathalpha}{letters}{"21}
\DeclareMathSymbol{\varepsilon}{\mathalpha}{letters}{"22}
\DeclareMathSymbol{\vartheta}{\mathalpha}{letters}{"23}
\DeclareMathSymbol{\varpi}{\mathalpha}{letters}{"24}
\DeclareMathSymbol{\varrho}{\mathalpha}{letters}{"25}
\DeclareMathSymbol{\varsigma}{\mathalpha}{letters}{"26}
\DeclareMathSymbol{\varphi}{\mathalpha}{letters}{"27}

\RedeclareMathOperator{\div}{\textbf{div}}

\renewcommand{\epsilon}{\varepsilon}

\newcommand{\half}{\ensuremath{\frac{1}{2}}}

\newcommand{\del}{\partial}

\newcommand{\DD}[2]{\frac{\partial #1}{\partial #2}}

\newcommand{\x}{\ensuremath{\vec x}\xspace}

\newcommand{\args}{\left(  t,\vec x\right)}

\newcommand\ie{i.e.\ }
\newcommand\eg{e.g.\ }

\let\origdoublepage\cleardoublepage
\renewcommand{\cleardoublepage}{%
	\clearpage{\pagestyle{empty}\origdoublepage}}

\makeatletter
\newcommand{\ChapterOutsidePart}{%
	\def\toclevel@chapter{-1}\def\toclevel@section{0}\def\toclevel@subsection{1}}
\newcommand{\ChapterInsidePart}{%
	\def\toclevel@chapter{0}\def\toclevel@section{1}\def\toclevel@subsection{2}}
\makeatother

\newcommand\wbref{target\xspace}

\newcommand{\dt}{\partial_t}
\newcommand{\tdt}{\frac{d}{dt}}
\newcommand{\dx}{\partial_x}
\newcommand{\dy}{\partial_y}

\newcommand{\q}{\vec q}
\newcommand{\s}{\vec s}
\newcommand{\f}{\vec f}
\newcommand{\g}{\vec g}

\newcommand{\Q}{\vec Q}

\newcommand{\F}{\vec F}

\begin{document}


\begin{frontmatter}

\title{High order well-balanced finite volume methods for multi-dimensional systems of hyperbolic balance laws
	}
\author[1]{Jonas P.\ Berberich\corref{cor1}}
\ead{jonas.berberich@mathematik.uni-wuerzburg.de}
\author[2]{Praveen Chandrashekar} 
\ead{praveen@tifrbng.res.in}
\author[1]{Christian Klingenberg}
\cortext[cor1]{Corresponding author: 
	Tel.: +49 931 31-88861;  
	Fax: +49 931 31-83494;}
\ead{klingen@mathematik.uni-wuerzburg.de}
\address[1]{Dept.~of Mathematics, Univ.~of W\"urzburg, Emil-Fischer-Stra{\ss}e 40, 97074 W\"urzburg, Germany}
\address[2]{TIFR Center for Applicable Mathematics, Bengaluru, Karnataka 560065, India}

\begin{abstract}
\reva{
We introduce a general framework for the construction of well-balanced finite volume methods for hyperbolic balance laws. We use the phrase \emph{well-balancing} in a broader sense, since our proposed method can be applied to exactly follow any solution of any system of hyperbolic balance laws in multiple spatial dimensions and not only time independent solutions. The solution has to be known a priori, either as an analytical expression or as discrete data. The proposed framework modifies the standard finite volume approach such that the well-balancing property is obtained and in case the method is high order accurate, this is maintained under our modification. We present numerical tests for the compressible Euler equations with and without gravity source term and with different equations of state, and for the equations of compressible ideal magnetohydrodynamics.
}
\end{abstract}


\end{frontmatter}


\section{Introduction}
\reva{Several problems in engineering and science are modeled by conservation principles and lead to non-linear partial differential equation which are hyperbolic. These equations can rarely be solved exactly and we must resort to some form of numerical approximation. One successful numerical approach for solving hyperbolic PDE is the finite volume method based on Godunov's idea \cite{Godunov1959}. Finite volume methods introduce discretization errors such that they are in general not exact on non-trivial solutions.} As soon as external forces enter the modeled system, a source term has to be added to the hyperbolic conservation laws turning these into hyperbolic balance laws.  \reva{\revb{Although numerical methods for hyperbolic balance laws might admit some discrete stationary states, they are in general grid dependent and different from the stationary states of the PDEs.} Standard numerical methods are not able to accurately maintain the latter solutions for long times, since they introduce discretization errors which are especially large on coarse grids.} This gives rise to the need to develop so-called well-balanced methods, i.e.\ methods which are designed to be exact on special stationary solutions of the system.

In the well-known shallow water equations with non-flat bottom topography, the most widely considered static state, which is the lake-at-rest solution, can be formulated in a closed form. This favors the construction of well-balanced methods for this system. There is a rich literature about well-balanced methods for shallow water equations (e.g.\ \cite{Brufau2002,Audusse2004,Bermudez1994,LeVeque1998b} and references therein) and related systems like the Ripa model (\cite{Desveaux2016b,Touma2015} and references therein). This includes high order methods for static \cite{Noelle2006} and non-static stationary states \cite{Noelle2007}. The relevance for methods for non-static stationary states has been pointed out in \cite{Xing2011}. For tsunami modeling applications, high order methods for shallow water equations on non-flat manifolds have been developed e.g.\ in \cite{Castro2018} considering the earth's surface geometry. For the Euler equations with gravitational potential, on the other hand, static solutions have to be found by solving a differential equation for density and pressure together with an equation of state (EoS). This makes the construction of well-balanced methods much more difficult and typically restricts the resulting method to some special cases. Many methods have been developed for some classes of hydrostatic states assuming an ideal gas EoS~ (\eg \cite{Cargo1994,LeVeque1999,LeVeque2011,Chandrashekar2015,Ghosh2016,Touma2016,Bermudez2017,Chertock2018a,Gaburro2018b} and references therein) and there are also high order methods, see e.g.\ \cite{Xing2013,Klingenberg2019,Gaburro2020b}.

While the hydrostatic equation for compressible Euler equations with gravity is basically one-dimensional, the spatial structure of this relation is much richer for compressible ideal magnetohydrodynamics (MHD) equations with gravity since it includes off-diagonal terms. In \cite{Fuchs2010}, a well-balanced method for MHD is derived to compute waves on the stationary background. This method is designed to balance isothermal hydrostatic states of the Euler equations together with a magnetic field, which satisfies certain stationarity conditions and is known a priori.
Part of this method, namely considering deviations to a background magnetic field, goes back to to Tanaka \cite{Tanaka1994} and is also used by Powell et al. \cite{Powell1999}. To do so, the background magnetic field is assumed to be static and free of rotation as well as divergence. 

There are different approaches to obtain the well-balanced property.
Some methods are based on a relaxation approach, in which the hydrostatic equation is included in the relaxation system \cite{Desveaux2014,Desveaux2016,Thomann2019}. Path-conservative methods are introduced in \cite{Pares2006,Castro2008,Gaburro2017}. Another widespread idea is the hydrostatic reconstruction, i.e.\ the reconstruction of variables which are constant if the system is in the considered stationary state. An early example of this method is \cite{Audusse2004} for the shallow water system; for Euler equations, this approach has \eg been used in \cite{Berberich2018,Berberich2019,Chandrashekar2015,Chertock2018a,Ghosh2016,Klingenberg2019,Touma2016}. The methods for Euler equations are restricted to a certain EoS and certain classes of hydrostatic solutions. For astrophysical applications, for example, this restriction is a severe limitation since the equations of state describing physics in the stellar interior are much more complex than the EoS of an ideal gas. 

More general methods have been developed in~\cite{Kaeppeli2014,Kaeppeli2016,Varma2019,Grosheintz2019}. The well-balanced methods introduced in these publications can be applied for any EoS. They are exact on certain hydrostatic solutions, and in all other cases they are exact for a second order approximation of the considered hydrostatic solution. In~\cite{Berberich2019}, a second order well-balanced method for Euler equations with gravity is introduced. This method can be applied for any hydrostatic solution of Euler equations with any EoS if the hydrostatic solution is known. The method is then exact up to machine precision and it has been extended to high order in~\cite{Klingenberg2019}. Notably, there is also an extension to stationary states with non-zero velocity in the same article. Similar techniques can be found in the context of numerical atmospheric modeling (e.g.\ \cite{Botta2004,Giraldo2008,Ghosh2016}). Those well-balanced schemes strongly rely on the structure of the discretized equations or the static solutions to be well-balanced.

The method we present in this paper is designed in the manner of the method in \cite{Berberich2019}. It uses the idea of hydrostatic reconstruction and a modification of the source term discretization to obtain the well-balanced property. The main point where it differs from all of the methods mentioned above is that our method is not restricted to a certain system of hyperbolic balance laws. Instead, we present a general framework which modifies finite volume methods for any hyperbolic conservation or balance laws such that they obtain the well-balancing property. Also, the method can be used to balance any known solution which is either given by an analytical expression or as discrete data. 
\revb{Additionally, unlike all the other well-balanced methods mentioned above which deal with time independent hydrostatic solutions, our well-balanced method can be used to follow known time-dependent solutions exactly}, as we will show in this paper. \reva{We will refer to the solution that is to be exactly captured by the scheme as the \emph{\wbref solution}. Depending on the application, the \wbref solution can be a time-independent hydrostatic solution, or some general time-dependent solution of interest.} Our method is also general in that it possible to combine it with other modules of a finite volume scheme: It can be applied on any grid system, with any numerical flux function, reconstruction routine, source term discretization, and ODE solver for time-discretization. It allows for high order in the sense that, if all these components are high order accurate, the resulting method is also high order accurate.

There are several applications in which the hydrostatic solution is known a priori, e.g., stellar astrophysics, but the EoS is often given in the form of a table. Consequently, hydrostatic solutions  can only be found numerically and are available in the form of discrete data. While methods which incorporate analytical expressions are not able to exactly maintain these hydrostatic solutions, it is very well possible with the methods in \cite{Berberich2019,Klingenberg2019} and the method we present in this paper. Especially, if we consider the better approximation of stellar structure which is given by a stationary state including rotation, our method can be applied to maintain this stationary solution. Another example from astrophysical application are rotating Keplerian disks. These are two-dimensional disks of matter which follows Newton's laws of motion in the gravitational field of a massive attractor. One way to describe this disk is a stationary solution of Euler equations with gravity including non-zero velocities. Since the velocity is not zero in such a solution, conventional well-balanced methods cannot preserve this solution. A special method designed for this application is presented in~\cite{Gaburro2018}. In this paper, we will show that our method is also able to preserve this solution on different grids. Besides the applicability to any system of hyperbolic balance laws, the balancing of moving and time-depending solutions is one of the key features of our method.

The rest of the paper is structured as follows. In \cref{sec:standardmethod3d} we introduce the standard finite volume framework  for systems of hyperbolic conservation laws in three spatial dimensions on arbitrary grids, \reva{but the numerical results will be presented only for 1-D and 2-D test cases}. In \cref{sec:method3d}, we introduce our general well-balanced modification for this framework. The well-balanced property we claim for our method is then shown in \cref{sec:proof}.
In \cref{sec:discrete_reference} the treatment of discrete \wbref solutions is discussed.
The validity of the well-balanced property also depends on a consistent choice of boundary conditions. Therefore, we add a discussion about well-balanced boundary conditions in \cref{sec:boundaries}. To emphasize how simple it is to add our method to an existing finite volume code, we comment on the implementation of the method in \cref{sec:implementation}. Finally, in \cref{sec:numerical_tests}, we show a variety of numerical tests. \reva{The range of applications goes from Euler equations  to ideal magnetohydrodynamics (MHD) equations.} They include classical well-balanced tests on the balance laws and also tests on the homogeneous hyperbolic conservation laws. Different equations of state are used for the Euler equations. We include a test in which the well-balanced solution is not analytically known but has been obtained numerically. Also, we present tests in which the well-balanced solution depends on time. We verify high order accuracy for solutions close to and far away from the well-balanced solution numerically. A simple example for using a \wbref solution which is obtained via numerical simulation is given.  The robustness of our approach is validated in a shock tube on a hydrostatic solution for Euler equations with gravity.
To show the efficiency of the method, we present CPU time comparisons of simulations with and without the well-balanced modification in \cref{sec:efficiency}.

\section{A standard finite volume method}
\label{sec:standardmethod3d}

\newcommand{\Fvec}{\mathcal{F}}

In this section we present the standard \reva{high} order finite volume framework for three-dimensional hyperbolic balance laws~\cite{Leveque1992,Toro2009}. Consider the 3-d system of hyperbolic balance laws
\begin{equation}
\dt\q(\x,t) + \nabla\cdot\Fvec(\q(\x,t))= \s(\q(\x,t),\x, t)
\label{eq:balance3d}
\end{equation}
with $\Fvec=(\f_1,\f_2,\f_3)$, where $\f_l$ is the flux in $l$-direction. The spatial domain is partitioned by a mesh consisting of $N$ non-overlapping control volumes.  For the $i$-th control volume $\Omega_i$ ($i\in\{1,\dots,N\}$), we define the cell-average
\begin{equation}
\Q_i(t) := \frac1{V_i}\int_{\Omega_i}\q(\x,t) d\x,
\end{equation}
where $V_i=|\Omega_i|$ is the control cell volume. Integrating \cref{eq:balance3d} over $\Omega_i$ and applying the divergence theorem yields an ordinary differential equation for $\Q_i$
\begin{equation}
\tdt\Q_i(t) + \frac1{V_i}\int_{\del\Omega_i}\Fvec(\q(\x,t))\cdot\vec n(\x)d\sigma = \frac1{V_i}\int_{\Omega_i}\s(\q(\x,t),\x)d\x,
\label{eq:standard_integrated_balance_law3d}
\end{equation}
\reva{which can be written as}
\begin{align}
\tdt\Q_i(t) = 
-\frac1{V_i}\sum_{k\in N(i)}\int_{\del\Omega_{ik}}\Fvec(\q(\x,t))\cdot\vec n(\x)d\sigma
+ \frac1{V_i}\int_{\Omega_i}\s(\q(\x,t),\x,t)d\x,
\label{eq:standard_perturbations_integrated_3d}
\end{align}
where $N(i)$ is the set of indices of all control volumes sharing an interface with $\Omega_i$ and $\partial\Omega_{ik}$ denotes the interface between  $\Omega_i$ and $\Omega_k$. For the discretization of the interface fluxes we use a \emph{numerical flux function}
$\F(\cdot,\cdot,\vec n)$
consistent with $\vec n\cdot\mathcal F$. The consistency conditions are Lipschitz continuity in the first two arguments and the relation $\F(\q,\q,\vec n)=\vec n\cdot\mathcal F(\q)$ for all unit vectors $\vec n$.
We apply this discretization to \cref{eq:standard_perturbations_integrated_3d} and obtain
\begin{align}
\tdt\Q_i(t) = -\frac1{V_i}\sum_{k\in N(i)}\int_{\del\Omega_{ik}}\F\left(\Q_i^\text{rec}(\x,t),\Q_k^\text{rec}(\x ,t),\vec n(\x)\right)d\sigma
+ \frac1{V_i}\int_{\Omega_i}\s\left(\Q_i^\text{rec}(\x,t),\x,t\right)d\x,
\label{eq:standard_integrated_3d_numflux}
\end{align} 
where the reconstructed functions $\Q_i^\text{rec},\Q_k^\text{rec}$ are obtained using a consistent conservative reconstruction routine on the cell average values $\Q$. Examples for popular consistent conservative reconstruction routines can be found in \cite{Harten1987,Liu1994,Levy2000,Toro2009}. In the next step we use numerical quadrature rules for the interface flux integral and a discretization of the source term integral. The semi-discrete method is then
\begin{align}
\tdt\Q_i(t) = -\frac1{V_i}\sum_{k\in N(i)}\left(\sum_{j=1}^{M}\omega_j\F\left(\Q_i^\text{rec}(\x_{ikj},t),\Q_k^\text{rec}(\x_{ikj} ,t),\vec n(\x_{ikj})\right)\right)
+ \frac1{ V_i}I_{\x\in\Omega_i}\left[\s\left(\Q_i^\text{rec},\x,t\right)\right].
\label{eq:standard_semi-discrete_3d}
\end{align}
$M$ is the number of quadrature points at the interfaces, $\x_{ikj}$ are the $M$ quadrature points at the $ik$ interface and $\omega_j$ are the corresponding weights. The symbol $I_{x\in\Omega}[\cdot]$ denotes a consistent discretization of the integral over the argument in the domain $\Omega$, 
\reva{
\[
I_{x\in\Omega}[\psi] \approx \int_\Omega \psi(\x) d\x
\]
}
The quadrature rules and source term discretizations we use in our tests are explained in the Appendix.

The semi-discrete scheme \cref{eq:standard_semi-discrete_3d} is $k$-th order accurate if the applied reconstruction routine, interface flux quadrature and source term discretization are all at least $k$-th order accurate.
It can then be evolved in time using a $k$-th order accurate ODE solver to obtain a $k$-th order accurate fully discrete scheme.

\section{The well-balanced modification of the standard finite volume method}
\label{sec:method3d}

In this section we will introduce a well-balanced modification for the three-dimensional finite volume method presented in \cref{sec:standardmethod3d}. Reducing it to one or two spatial dimensions is straight forward.

Let $\tilde{\q}$ be a given continuous and sufficiently smooth solution of \cref{eq:balance3d}. Plugging this \wbref solution $\tilde\q$ into \cref{eq:standard_integrated_balance_law3d} we get
\begin{equation}
\tdt\tilde\Q_i(t) + \frac1{V_i}\int_{\del\Omega_i}\Fvec(\tilde\q(\x,t))\cdot\vec n(\x)d\sigma = \frac1{V_i}\int_{\Omega_i}\s(\tilde\q(\x,t),\x)d\x,
\label{eq:integrated_balance_law3d_stationary}
\end{equation}
where $\tilde{\Q}_i$ is the average of the \wbref solution $\tilde \q$ in the $i$-th control volume.
In the next step we subtract \cref{eq:integrated_balance_law3d_stationary} from \cref{eq:standard_integrated_balance_law3d} to obtain
\begin{equation}
\tdt\Q_i(t)-\tdt\tilde\Q_i(t) + \frac1{V_i}\int_{\del\Omega_i}(\Fvec(\q(\x,t)) 
-\Fvec(\tilde\q(\x,t)))\cdot\vec n(\x)d\sigma
= \frac1{V_i}\int_{\Omega_i}\s(\q(\x,t),\x,t)-\s(\tilde\q(\x,t),\x,t)d\x.
\label{eq:integrated_balance_law3d_substracted}
\end{equation}
Now, let us rewrite \cref{eq:integrated_balance_law3d_substracted} in terms of the deviation from the \wbref solution
\begin{align}
\Delta \q &:= \q-\tilde{\q},& \Delta \Q &:= \Q-\tilde{\Q}. 
\end{align}
This yields
\begin{align}
\tdt(\Delta\Q_i(t)) = 
& 
-\frac1{V_i}\sum_{k\in N(i)}\int_{\del\Omega_{ik}}(\Fvec((\Delta\q+\tilde{\q})(\x,t)) 
-\Fvec(\tilde\q(\x,t)))\cdot\vec n(\x)d\sigma\nonumber\\
& 
+ \frac1{V_i}\int_{\Omega_i}\s((\Delta\q+\tilde{\q})(\x,t),\x,t)-\s(\tilde\q(\x,t),\x,t)d\x,
\label{eq:perturbations_integrated_3d}
\end{align}
where $N(i)$ is the set of indices of all control volumes sharing an interface with $\Omega_i$.
At this point, we start to discretize. For that we define a numerical flux difference approximation
\begin{equation}
\Delta\hat\Fvec\left(\Delta\Q^L,\Delta\Q^R,\tilde\q,\vec n\right) := 
\F(\Delta\Q^L+\tilde\q,\Delta\Q^R+\tilde\q,\vec n) - \vec n\cdot\mathcal F(\tilde{\q})
\approx \vec n\cdot(\Fvec(\Delta\q+\tilde{\q})-\Fvec(\tilde{\q})),
\label{eq:flux_difference_3d}
\end{equation}
where $\F(\cdot,\cdot,\vec n)$ is a numerical flux function consistent with $\vec n\cdot\mathcal F$. 
We apply this discretization to \cref{eq:perturbations_integrated_3d} and obtain
\begin{align}
\tdt(\Delta\Q_i(t)) = &-\frac1{V_i}\sum_{k\in N(i)}\int_{\del\Omega_{ik}}\Delta\hat\Fvec\left(\Delta\Q_i^\text{rec}(\x,t),\Delta\Q_k^\text{rec}(\x ,t),\tilde\q(\x ,t),\vec n(\x)\right)d\sigma\nonumber\\
&+ \frac1{V_i}\int_{\Omega_i}\s((\Delta\Q_i^\text{rec}+\tilde{\q})(\x,t),\x,t)-\s(\tilde\q(\x,t),\x,t)d\x,
\label{eq:perturbations_integrated_3d_numflux}
\end{align} 
where the reconstructed functions $\Delta\Q_i^\text{rec},\Delta\Q_k^\text{rec}$ are obtained using a consistent conservative reconstruction routine on the cell average values $\Delta\Q$. 
\reva{Note that using a consistent reconstruction on the deviations $\Delta\Q$ is equivalent to the application of a hydrostatic reconstruction as \eg in \cite{Audusse2004,Berberich2020c}. A typical hydrostatic reconstruction consists of a transformation to hydrostatic variables (\ie a set of variables that is constant in the hydrostatic case), a consistent reconstruction, and a transformation back to conservative variables. In our description of the method the deviations $\Delta \Q$ correspond to the hydrostatic states. Hence, no transformation is required.}
In the next step we use numerical quadrature rules for the interface flux integral and a discretization of the source term integral. The semi-discrete method is then
\begin{align}
\tdt(\Delta\Q_i(t)) = &-\frac1{V_i}\sum_{k\in N(i)}\left(\sum_{j=1}^{M}\omega_j\Delta\hat\Fvec\left(\Delta\Q_i^\text{rec}(\x_{ikj},t),\Delta\Q_k^\text{rec}(\x_{ikj} ,t),\tilde\q(\x_{ikj} ,t),\vec n(\x_{ikj})\right)\right)\nonumber\\
&+ \frac1{ V_i}I_{\x\in\Omega_i}\left[\s((\Delta\Q_i^\text{rec}+\tilde\q)(\x,t),\x,t)\right]
-\frac1{ V_i}I_{\x\in\Omega_i}\left[\s(\tilde\q(\x,t),\x,t)\right].
\label{eq:semi-discrete_3d}
\end{align}
where the notations used are as in \cref{eq:standard_semi-discrete_3d}.
\revb{
	If the source term $\s$ is linear in the first argument, we can use the following relation
	\begin{equation}
	\s((\tilde{\q}+\Delta\q)\args,\x,t)-\s(\tilde{\q}\args,\x,t) = \s(\Delta\q\args,\x,t).
	\end{equation}
	Due to the linearity of the corresponding source term discretizations, this relation then also holds for the discretized source terms, which further simplifies the scheme.
	For example, this is the case for the gravitational source term in Euler or ideal MHD equations and the bottom topography source term in the shallow water equations.
}

As in the standard method, this semi-discrete scheme \cref{eq:semi-discrete_3d} is $k$-th order accurate if the applied reconstruction routine, interface flux quadrature and source term discretization are all at least $k$-th order accurate.
It can then be evolved in time using \revb{an at least} $k$-th order accurate ODE solver to obtain a $k$-th order accurate fully discrete scheme.

\revb{
\begin{myremark}
    \label{rem:discontinuous}
    In the description of the method, we assume the \wbref solution $\tilde \q$ to be smooth. In the case of discontinuous $\tilde \q$, the two values $\tilde\q^L$ and $\tilde\q^R$ which are different if a discontinuity is present at the interface, have to given instead of only one value. \Cref{eq:flux_difference_3d} has then to be modified to
    \begin{equation}
        \Delta\hat\Fvec\left(\Delta\Q^L,\Delta\Q^R,\tilde\q^L,\tilde\q^R,\vec n\right) :=
        \F(\Delta\Q^L+\tilde\q^L,\Delta\Q^R+\tilde\q^R,\vec n) - \F(\tilde\q^L,\tilde\q^R,\vec n),
        \label{eq:flux_difference_3d_discontinuous}
    \end{equation}
    \ie the numerical flux function is also applied to the \wbref solution. If $\tilde\q$ is continuous, \cref{eq:flux_difference_3d_discontinuous} reduces to \cref{eq:flux_difference_3d} due to the consistency of the numerical flux function $\F$. However, in the case of a discontinuous \wbref solution $\tilde\q$ no high order convergence can be expected and the computational cost of the well-balanced modification increases.
\end{myremark}
}

\section{Proof of the well-balanced property}
\label{sec:proof}

In this section we show the well-balanced property of our method.
\begin{mytheorem}
	\label{thm:wbproperty}
	The modified finite volume method introduced in \cref{sec:method3d}
	satisfies the following property: If
	\begin{equation}
		\Delta\Q_{i}=0\qquad\forall i\in \{1,\dots,N\}
	\end{equation} 
	at initial time, then this holds for all $t>0$. Consequently, if the initial condition $\Q_i(t=0)$, $i=1,\dots,N$, equals the cell averages of the \wbref solution  $\tilde\Q_i(t=0)$, $i=1,\dots,N$, the computed solution equals the \wbref solution for all time.
\end{mytheorem}
\emph{Proof:} 
Let $\Delta \Q_i=0$ for all $i\in \{1,\dots,N\}$. The consistency of the applied reconstruction leads to 
$\Delta \Q_i^\text{rec}\equiv 0$
at all flux quadrature points. The flux consistency then yields 
\begin{equation}
\label{eq:fluxconsistency}
\Delta\hat\Fvec\left(\Delta\Q^L,\Delta\Q^R,\tilde\q,\vec n\right)= 
\Delta\hat\Fvec\left(0,0,\tilde\q,\vec n\right)= 
\F(\tilde\q,\tilde\q,\vec n) - \vec n\cdot\mathcal F(\tilde{\q})
=\vec n\cdot\mathcal F(\tilde{\q}) - \vec n\cdot\mathcal F(\tilde{\q})
=0.
\end{equation}
Now, consider the contribution from the source term: With $\Delta\Q_i=0$ the source term discretization in \cref{eq:semi-discrete_3d} reduces to
\begin{equation}
I_{\x\in\Omega_i}\left[\s((\Delta\Q_i^\text{rec}+\tilde{\q})(\x,t),\x,t)\right]
-I_{\x\in\Omega_i}\left[\s(\tilde\q(\x,t),\x,t)\right]
=I_{\x\in\Omega_i}\left[\s(\tilde{\q}(\x,t),\x,t)\right]
-I_{\x\in\Omega_i}\left[\s(\tilde\q(\x,t),\x,t)\right]=0.
\end{equation}
We have shown that the right hand side in \cref{eq:semi-discrete_3d} vanishes and thus the initial data $\Delta\Q_i=0$ are conserved for all time. The second part of the theorem follows easily.
\hfill\qed

\revb{
\begin{myremark}
	\label{rem:evolution}
	If a stationary solution is chosen as \wbref solution (which is the case for classical well-balancing applications), the time derivative of the \wbref solution vanishes by definition. This leads to $\tdt\Q_i=\tdt(\Delta\Q_i)$. The described method can then also be used to directly evolve the $\Q_i$ in time instead of $\Delta\Q_i$.
\end{myremark}
}

\section{The treatment of discrete \wbref solutions}
\label{sec:discrete_reference}

\reva{
The \wbref solution, introduced in \cref{sec:method3d} for the well-balanced modification, has not to be known analytically, it can also be given only in the form of discrete data. 
}
For consistency with the system of balance laws \cref{eq:balance3d} it is important that the discrete data which are used for the \wbref solution converge to a solution of \cref{eq:balance3d} when the computational grid is refined. To ensure the high order of our method, this convergence should also be of high order. When only discrete data are given we need a method to compute values at the interface quadrature points and cell-averaged values for the grid on which we use our well-balanced method. \revb{Note that the \wbref solution obtained from the discrete data are in general only approximate solutions of the hyperbolic system, but we assume that it is a sufficiently accurate approximation to the true solution of the PDE. The well-balanced modification introduced in this article is then exact \emph{on this approximate solution}.}

 \revb{
 Depending on the form in which the discrete \wbref data are given, \eg point values or cell averages, and the required order of accuracy, different methods for this reconstruction have to be used. Also, in the case of a discontinuous \wbref solution
 (see \cref{rem:discontinuous}), two interface values of $\tilde\q$ have to be determined for each interface quadrature point instead of only one.
In the following we give two examples how cell average and interface values of a \wbref solution can be obtained from discrete data.}

\subsection{Example 1: pointwise 1-D data on a fine grid}
\label{sec:discrete_reference_example1}

\reva{An important application of our method is in  computing hydrostatic solutions of Euler equations and perturbations around such solutions.} Especially in physical applications with complex EoS hydrostatic solutions have to be obtained by numerical methods.
Even for multi-dimensional simulations, the underlying hydrostatic solution can be one-dimensional in its nature (see cases (a) and (b) below).
Now assume such a numerically approximated hydrostatic solution in one spatial dimension is given in the form of point values $\q^{hs}_i$, $i=1,2,\dots,L$ (in conservative variables) on a fine equidistant grid. Assume it is supposed to be used in a two-dimensional third order accurate modified finite volume method as introduced in \cref{sec:method3d} on a Cartesian grid. For that, in a first step, we use a cubic 
spline interpolation to construct a continuous function 
$\q^{hs}_{1-D}(x)$ (e.g.\ \cite{Press1992}). This function is then extended to two spatial dimensions. How this is done depends on the symmetry of the 2-D problem which allowed the reduction to a 1-D hydrostatic solution. We consider two different cases:
\begin{itemize}
	\item[(a)] Suppose we have an essentially 1-D hydrostatic solution where the gravitational force is at an angle $\alpha\in[0,2\pi)$ to the $x$-axis.
	Then we extend the one-dimensional hydrostatic solution $\q^{hs}_{1-D}$ to a two-dimensional solution via $\q^{hs}(\x):=\q^{hs}_{1-D}(x\cos(\alpha)+y\sin(\alpha))$.
	\item[(b)] Suppose we have a radial hydrostatic solution with a gravity vector pointing towards the center $\x_\text{center}$. In that case we extend the hydrostatic solution back to 2-D by setting $\q^{hs}(\x):=\q^{hs}_{1-D}(\|\x-\x_\text{center}\|)$.
\end{itemize}
The values for the \wbref solution at interface quadrature points can then be evaluated pointwise as $\tilde\q(\x,t)=\q^{hs}(\x)$. The cell average values of the \wbref solution are computed using a third order accurate 2-D Gau\ss--Legendre quadrature rule. Method (a) is used in the test in \cref{sec:integrated_hystat}.
\subsection{Example 2: \wbref solution from a highly resolved finite volume simulation}
\label{sec:discrete_reference_example2}
 
Consider a \wbref solution given as numerical solution of a finite volume simulation on a two-dimensional structured static grid with curvilinear coordinates as described in \ref{appendix:curvilinear}. In this case the data are given as cell averages $\hat{\Q}^n_{ij}$, instead of point values. Assume this \wbref solution is supposed to be used in a $k$-th order accurate two-dimensional modified finite volume method as described in \cref{sec:method3d}. The grid to be used is a coarser version of the grid on which the \wbref solution has been computed, such that all interfaces on the coarse grid coincide with interfaces of the fine grid. For each time step $t^n$ of the stored \wbref data we map the fine grid on the coarse grid with 
\begin{equation}
\tilde\Q_{ij}(t^n)=\frac1{\tilde V_{ij}}\sum_{\hat\x_{kl}\in\tilde\Omega_{ij}}\hat V_{kl}\hat{\Q}^n_{kl},
\end{equation}
where all quantities with $\tilde{\cdot}$ correspond to the coarse grid and all quantities with $\hat{\cdot}$ correspond to the fine grid. The values of $\tilde\Q_{ij}$ at intermediate times are obtained via a $k$-th order accurate interpolation in time.

The value of the \wbref solution at all quadrature points required in the scheme are obtained using a $k$-th order accurate interpolation on the cell-centered point values $\tilde\q^\text{rec}(\x_{ij})$. Those cell-centered values are obtained using a $k$-th order accurate conservative reconstruction on the cell-averages of the cell-average values $\tilde{\Q}_{ij}$.
This method is applied in a numerical test in \cref{sec:mhd_vortex_numref} for $k=3$.

\section{Boundary conditions}
\label{sec:boundaries}
In the previous sections (including the proof of the well-balanced property) we omitted to include boundary conditions in the discussion. Yet, the validity of the well-balanced property also depends on the correct choice of boundary conditions. \revb{In this section, we will describe some boundary conditions which are compatible with the well-balancing property and support the potentially high order accuracy of the scheme.} Some of the proposed numerical boundary conditions require knowledge of the \wbref solution outside the domain. \revb{If this is not available, one can simply extrapolate the \wbref solution to the ghost cells using a sufficiently high order accurate extrapolation. The well-balanced scheme will in that case exactly balance the \wbref solution including the extrapolated values. Consequently, this will not affect the well-balanced property nor order of accuracy.}

\noindent\emph{Extrapolation boundary conditions}: One way to treat boundaries is the extrapolation of data in the domain to ghost cells. This can be done with high order to support the high order of the applied scheme. In our method we extrapolate the deviations $\Delta\Q$. In the case that $\Delta\Q_i=0$ holds for all control volumes in the domain, this will also be true for the extrapolated states. Hence, the well-balanced property also holds at the boundary. Extrapolation boundary conditions for one and two spatial dimensions with different orders of accuracy can be found in \ref{appendix:boundaries}.

\revb{
\noindent\emph{Wall and periodic boundary conditions}
Periodic boundary conditions and wall boundary conditions as \eg described in \cite{Blazek2015} can be applied to the deviations $\Delta\Q$. In the case of $\Delta \Q=0$ for each cell in the domain, this implies $\Delta\Q=0$ for each ghost cell. Hence, these boundary conditions are consistent with the well-balanced property.
}

\section{Notes on the implementation}
\label{sec:implementation}

We have seen that our well-balanced method can be applied for a wide range of problems. In this section we will discuss another property useful \reva{for applications: The fact that the method is also easy to be implemented}. This holds especially if there is an existing finite volume code \reva{evolving \cref{eq:standard_semi-discrete_3d}. It can be easily modified to evolve \cref{eq:semi-discrete_3d} and hence obtain a well-balancing capability. The changes that should be introduced in a typical finite volume scheme are the following:}
\begin{enumerate}
	\item Implement a function \reva{able to return the \wbref solution $\tilde \q$ in any point $(\x,t)$, or a structure containing the cell averages and all the values needed at each quadrature point} if $\tilde\q$ is time-independent.
	\item \reva{When the main routine receives the initial cell-averages $\Q_i$, it should transform it to the deviations $\Delta \Q_i=\Q_i-\tilde{\Q_i}$ and work directly with the deviation.} 
	\item In the routine evaluating the numerical flux, \reva{it has to be computed} $\Delta\hat\Fvec\left(\Delta\Q^L,\Delta\Q^R,\tilde\q,\vec n\right)$ instead of $\vec n\cdot\Fvec(\Delta\q+\tilde{\q})$. Basically, this just means subtracting the exact flux after \reva{the} evaluation of the \reva{standard} numerical flux.
	\item This step is only necessary if the source term is not linear in $\q$: Evaluate the source term at the states $\tilde\Q_i+\Delta\Q_i$ and $\tilde{\Q_i}$. \reva{The difference of these source terms is computed and used to evolve the approximate solution as described in \cref{eq:semi-discrete_3d}.}
\end{enumerate}
Let us remind of \cref{rem:evolution} and point out that an alternative implementation could also evolve $\Q$ instead of $\Delta\Q$ in time if $\tilde{\q}$ is time-independent. \reva{In the computer implementation, this can be the easier and computationally efficient approach than evolving the deviations.}
\section{Numerical tests of the scheme}
\label{sec:numerical_tests}

\reva{
\subsection{Hyperbolic systems used in the tests}
Since the well-balancing procedure introduced in this article is applicable to any hyperbolic balance law, we present numerical experiments for two different hyperbolic systems. They are introduced in the following.
}
\subsubsection{Compressible Euler equations with gravity source term}
The 2-D compressible Euler equations which model the balance laws of mass, momentum, and energy under the influence of gravity are given by
\begin{equation}
\dt{\q}+ \dx{\f} + \dy{\g} = \s,
\label{eq:eul2d_cartesian}
\end{equation}
where the conserved variables, fluxes and source terms are
\begin{align}
\label{eq:euler2d}
\q &= \begin{bmatrix}
\rho \\ \rho u \\ \rho v \\ E \end{bmatrix}, &
\f &= \begin{bmatrix}
\rho u \\
p + \rho u^2 \\
\rho u v \\
(E+p)u \end{bmatrix}, & 
\g &= \begin{bmatrix}
\rho v \\
\rho u v \\
p + \rho v^2 \\
(E+p)v \end{bmatrix}, &
\s &= \begin{bmatrix}
0 \\
-\rho \dx{\phi}\\
-\rho \dy{\phi}\\
0
\end{bmatrix}
\end{align}
with $\rho,p>0$.
Moreover, $E=\rho\varepsilon + \tfrac{1}{2}\rho|\vec v|^2 + \rho\phi$ is the total energy per unit volume with the velocity $\vec v=(u,v)^T$ and specific internal energy $\epsilon$. The scalar function $\phi$ is a given gravitational potential. An additional relation between density, pressure, and specific internal energy is given in the form of an equation of state (EoS). In our tests we will use the ideal gas EoS
\begin{equation}
p=(\gamma-1)\rho\epsilon
\end{equation} 
with $\gamma=1.4$, although our well-balanced method can be applied for Euler equations with any EoS.

The 2-D Euler equations can be reduced to 1-D Euler equations by setting $\g=0$ and removing the $\rho v$ equation. It can be reduced to homogeneous Euler equations by setting $\s=0$.

\subsubsection{Homogeneous compressible ideal magnetohydrodynamics}

The 2-D compressible ideal magnetohydrodynamics (MHD) equations which model the conservation of mass, momentum, magnetic field, and energy are given by
\begin{equation}
\dt{\q}+ \dx{\f} + \dy{\g} = 0.
\label{eq:mhd2d_cartesian}
\end{equation}
The conserved variables and fluxes are
\begin{align}
\label{eq:mhd2d}
\q &= \begin{bmatrix}
\rho \\ \rho u \\ \rho v \\ B_x \\ B_y \\ E \end{bmatrix}, &
\f &= \begin{bmatrix}
\rho u \\
\rho u^2 + p + \frac12 (B_y^2-B_x^2) \\
\rho u v -B_xB_y\\
0\\
B_yu-vB_x\\
u(E+p+\frac12B_y^2-\frac12B_x^2)-vB_xB_y \end{bmatrix},&
\g &= \begin{bmatrix}
\rho v \\
\rho u v - B_x B_y\\
\rho v^2 + p + \frac12 (B_x^2-B_y^2) \\
B_xv-uB_y\\
0\\
v(E+p+\frac12B_x^2-\frac12B_y^2)-uB_xB_y \end{bmatrix},
\end{align}
where $B_x$, $B_y$ are the $x$- and $y$-component of the magnetic field. The total energy is $E=\rho\varepsilon + \tfrac{1}{2}\rho|\vec v|^2 + \frac12(B_x^2+B_y^2)$. All other quantities are defined as for the Euler equations.  We use the same EoS as for the Euler equations.

One can also define 2-D compressible ideal MHD equations such that they include the $\rho w$ and $B_z$ components. This is in principle reasonable due to the genuine three-dimensional interactions between velocity and magnetic field. In our tests we set $\rho w$ and $B_z$ to zero and there is no difference if we omit the corresponding equations.

\subsection{Code and numerical methods}
\revb{We test the methods proposed in this paper using a finite volume code implemented in Python. The code is \reva{built} in a modular way, such that different \reva{schemes} can be easily implemented.} \reva{For brevity, we only give a short description of the methods in this section and the interested reader can find details in the \ref{appendix:detail}. 
Note that our well-balanced method is not restricted to the methods we choose to use in the tests and one can for example also use other reconstruction methods or quadrature formulae. One can also apply numerical flux functions designed for special problems (e.g.\ a low Mach number compliant method for Euler equations like in \cite{Barsukow2017,Berberich2020b}).} \\
\noindent\emph{Grids:} \reva{
The domain in two-dimensional problems is discretized using a structured grid and in some tests, we use curvilinear grids. The implementation of these grids restricts the overall method to only second order accuracy. Note that this is not a general statement about the proposed method, but this restriction is due to the special implementation of the grid in our code.  Higher order accuracy can only be achieved with a Cartesian grid in our code.} \\
\noindent\emph{Numerical flux function:} As numerical flux function we use the local Lax--Friedrichs flux (e.g.\ \cite{Leveque1992}), since it is simple and can be applied for any hyperbolic system. In some tests, we use the Roe's approximate Riemann solver for Euler equations~\cite{Roe1981} to obtain more accurate results. \\
\noindent\emph{First order method:} To formally obtain a first order method, we use constant reconstruction to obtain the interface values. \reva{The numerical fluxes are computed at the center of the interfaces and the source term is evaluated at the cell-center. For the gravity source term used in our tests, we need the gradient $\nabla \phi$ of the given gravitational potential which is evaluated at the cell center using analytical differentiation.} \\
\noindent\emph{Second order method:} To formally obtain a second order method, we use a conservative linear reconstruction (e.g.\ \cite{Toro2009}) with a minmod limiter (e.g.\ \cite{Toro2009}) to obtain the interface values. This is the only difference to the first order method. \\
\noindent\emph{Third order method:} 
To formally obtain a third order method, we use a conservative CWENO3 (\cite{Levy2000} for 1-D, \cite{Levy2000b} for 2-D) reconstruction to obtain the interface values. In the two-dimensional case, the numerical fluxes are evaluated at the Gau\ss--Legendre quadrature points. 
To compute the source term we multiply the CWENO3 polynomials in momentum with the interpolation polynomial of $\del_x \phi$ or $\del_y\phi$ respectively. \reva{This results in a polynomial function of the spatial variables whose cell average can be computed by analytical integration.}
\\\noindent\emph{Seventh order method:} 
To formally obtain a one-dimensional seventh order method, we use a conservative CWENO7 reconstruction \cite{Cravero2018} to obtain the interface values. To compute the source term we multiply the CWENO7 polynomial in momentum with the interpolation polynomial of $\del_x \phi$. The resulting source term polynomial is cell-averaged using analytical integration. \\
\noindent\emph{Boundary conditions:} If the setup has periodic character we use periodic boundary conditions. Otherwise we extrapolate the states to ghost cells with a sufficiently high spatial order. If we use the third order method, for example, we extrapolate using parabolas. \\
\revb{
\noindent\emph{Time-stepping:} The first order accurate scheme is evolved in time using the explicit forward Euler method. For the second and third order accurate semi-discrete scheme we use the explicit third order, four stage Runge--Kutta method from \cite{Kraaijevanger1991} and the explicit tenth order, 17 stage Runge--Kutta method from \cite{Feagin2007} is used in the seventh order accurate method.
}


\subsection{1-D isothermal hydrostatic solution of the Euler system with gravity}
\label{sec:1d_isothermal}
We consider an isothermal hydrostatic solution of the 1-D compressible \revb{Euler equations with the ideal gas law and the gravitational source term given by}
\begin{align}
\phi(x)&=\sin(2\pi x),& 
\tilde\rho(x)&=\tilde p(x)=\exp(-\phi(x)),& 
\tilde u&\equiv0.
\label{eq:1d_isothermal}
\end{align}
We set these data on a 1-D grid with 128 grid cells on the domain $[0,1]$. These initial data are evolved up to the final time $t=2$ using the first, second, and third order method with the standard method and the well-balanced method each. In the well-balanced method we set the initial data \cref{eq:1d_isothermal} as time-independent \wbref solution. The $L^1$-errors at final time compared to the initial grid can be seen in \cref{tab:1d_isothermal}. We see that there is no error when the well-balanced method is applied. 
\begin{table}
	\centering
	\caption{
		\label{tab:1d_isothermal}
		$L^1$-errors for an isothermal hydrostatic solution of the Euler equations after time $t=2$ computed using the standard (Std) and well-balanced (WB) methods with different orders of accuracy (O$m$ for $m$-th order). The setup is described in \cref{sec:1d_isothermal}.
	}
	\setlength\tabcolsep{5pt}
	\begin{tabular}{| c | c | c | c | c | c | c | c | c |}
		\hline
		error in & Std-O1 & WB-O1 & Std-O2 & WB-O2 & Std-O3 & WB-O3 & Std-O7 & WB-O7   \\
		\hline
		$\rho$  & 1.19e-01 & 0.00e+00 & 4.60e-04 & 0.00e+00 & 9.72e-05 & 0.00e+00 & 1.28e-09 & 0.00e+00  
		\\
		$\rho u$& 2.18e-02 & 0.00e+00 & 7.24e-04 & 0.00e+00 & 1.50e-04 & 0.00e+00 & 1.40e-09 & 0.00e+00  
		\\
		$E$     & 1.64e-01 & 0.00e+00 & 3.20e-03 & 0.00e+00 & 3.92e-04 & 0.00e+00 & 3.99e-09 & 0.00e+00 
		\\
		\hline
	\end{tabular}
	\setlength\tabcolsep{6pt}
\end{table}
\begin{myremark}
	Most other well-balanced methods balance fluxes against the source term which leads to machine errors.
	In our method we balance fluxes against fluxes and source term against source term. Thus, the differences can cancel out exactly and the error can be exactly zero.
\end{myremark}
\subsection{Perturbed 1-D isothermal hydrostatic solution of the Euler system with gravity}
\label{sec:1d_perturbation}
\begin{figure}[tbph]
	\centering
	\includegraphics[scale=1]{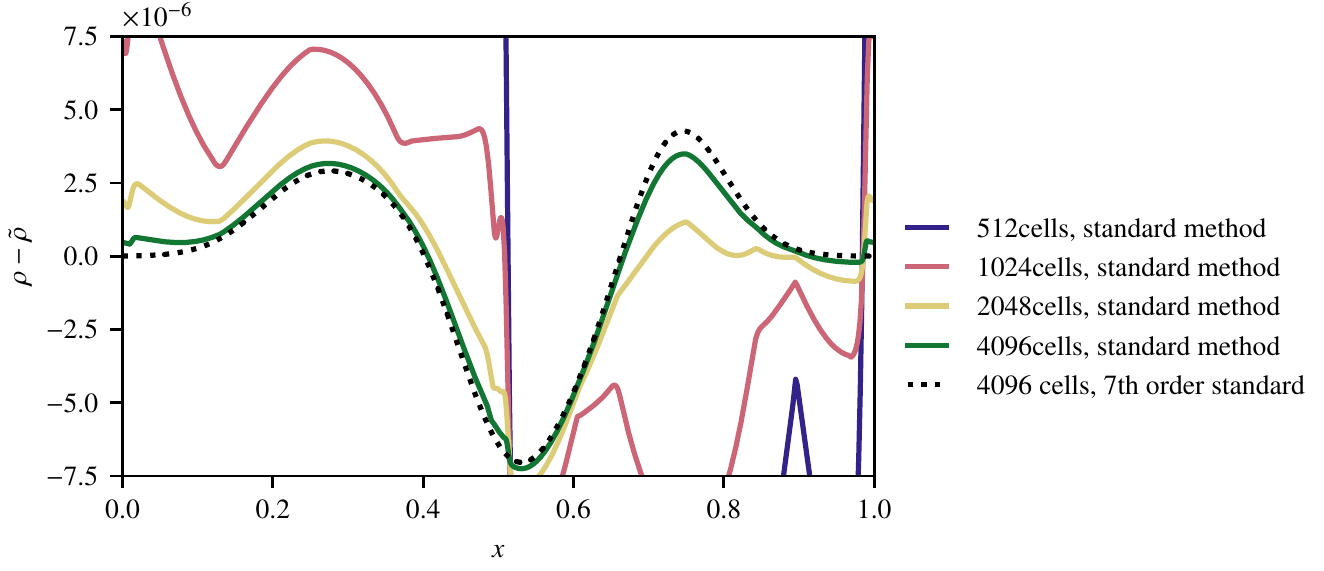}
	\caption{Perturbation on a hydrostatic atmosphere. The test setup is described in \cref{sec:1d_perturbation}. The density deviation from the hydrostatic background is shown at time $t=0.2$ for the perturbation $\eta=1e-5$. The second order (if not stated explicitly) standard method is used with different resolutions.}
	\label{fig:isothermal_1e-5_nowb}
\end{figure}

\begin{figure}[tbph]
	\centering
	\includegraphics[scale=1]{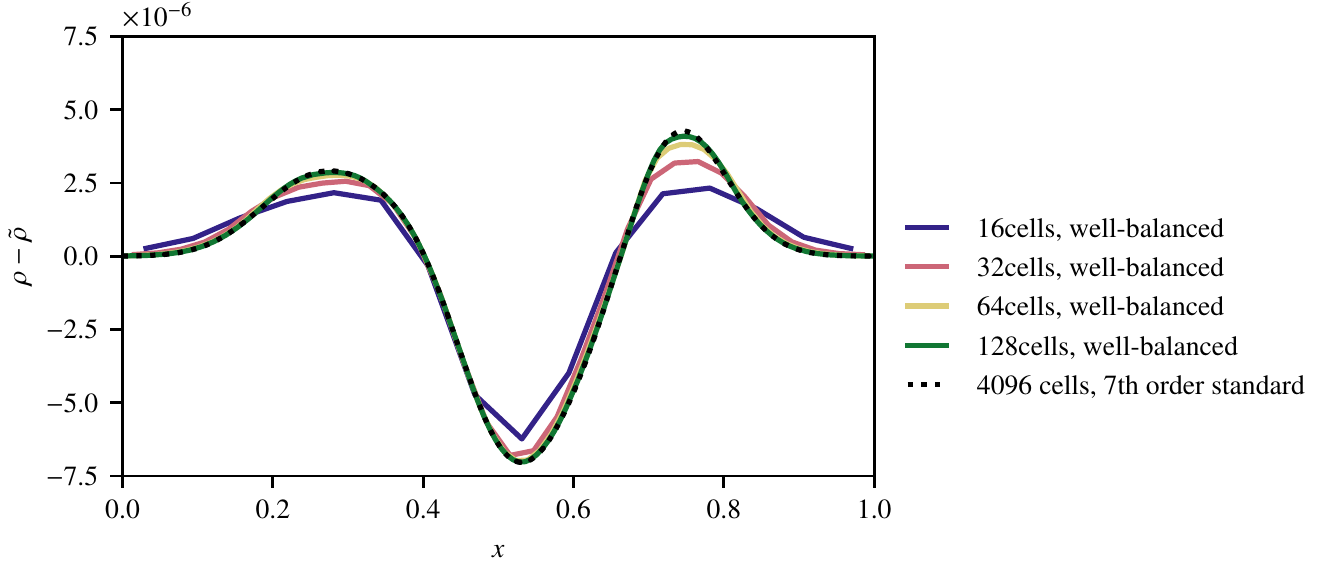}
	\caption{Perturbation on a hydrostatic atmosphere. The test setup is described in \cref{sec:1d_perturbation}. The density deviation from the hydrostatic background is shown at time $t=0.2$ for the perturbation $\eta=1e-5$. The second order well-balanced method is used with different resolutions (if not stated explicitly).}
	\label{fig:isothermal_1e-5_wb}
\end{figure}

We add a perturbation to the pressure such that our initial conditions are
\begin{align}
\rho(x)&=\tilde{\rho}(x),& 
u(x)   &=\tilde u(x), & 
p(x)   &=\tilde{p}(x)+\eta \exp\left(-100\left(x-\half\right)^2\right)
\end{align}
in the domain $[0,1]$.
We choose $\eta=0.1$ to test the convergence of our method. 
We evolve this initial setup up to time $t=0.2$ using our well-balanced method (first to third order and seventh order). 
The results and convergence rates are shown in \cref{tab:1d_perturbation}. As a reference solution we use a numerical solution computed with the seventh order standard scheme on a grid with $4096$ cells. All convergence rates match our expectations. The convergence rate for the seventh order scheme drops in the last step, since the error approaches machine precision. In \cref{fig:isothermal_1e-5_nowb,fig:isothermal_1e-5_wb}, density deviations at time $t=0.2$ for the test with $\eta=10^{-5}$ are shown. In \cref{fig:isothermal_1e-5_nowb} we see, that the discretization error on the hydrostatic background dominates the total error when the second order standard method is used with a low resolution. For higher resolutions, the perturbation can be resolved correctly. In \cref{fig:isothermal_1e-5_wb} on the other hand, it gets evident that the second order well-balanced method is capable of correctly resolving the perturbation on a coarse grid. This is due to the fact that there are no discretization errors on the hydrostatic background as we have already seen in \cref{sec:1d_isothermal}.


\begin{table}
	\centering	
	\caption{
		\label{tab:1d_perturbation}
		$L^1$-errors and convergence rates in total energy for a small pressure perturbation on an isothermal hydrostatic solution of the Euler equations after time $t=0.2$. Different well-balanced methods are used. The setup is described in \cref{sec:1d_perturbation}.
	}
	
	\begin{tabular}{| c | c c | c c | c c | c c|}
	    \hline
		\multirow{2}{*}{$N$} & \multicolumn{2}{c|}{WB-O1} & \multicolumn{2}{c|}{WB-O2} & \multicolumn{2}{c|}{WB-O3} & \multicolumn{2}{c|}{WB-O7}   \\
		 & $E$ error & rate& $E$ error & rate& $E$ error & rate& $E$ error & rate \\
		\hline
		256  & 5.73e-03 & \multirow{2}{*}{0.9} & 5.98e-05 & \multirow{2}{*}{2.0} & 6.93e-05 & \multirow{2}{*}{2.6} & 1.32e-09 & \multirow{2}{*}{6.5}\\
		512  & 3.08e-03 & \multirow{2}{*}{0.9} & 1.49e-05 & \multirow{2}{*}{2.0} & 1.18e-05 & \multirow{2}{*}{2.7} & 1.46e-11 &  \multirow{2}{*}{6.6}\\
		1024 & 1.60e-03 & \multirow{2}{*}{1.0} & 3.73e-06 & \multirow{2}{*}{2.0} & 1.78e-06 & \multirow{2}{*}{2.8} & 1.46e-13 &  \multirow{2}{*}{4.1}\\
		2048 & 8.15e-04 &                      & 9.36e-07 &                      & 2.53e-07 &                      & 8.41e-15 & \\
	    \hline
	\end{tabular}

\end{table}
\subsection{Riemann problem on a 1-D isothermal hydrostatic solution of the Euler system with gravity}
\label{sec:robustness}
\begin{figure}[h!]
	\centering
	\includegraphics[scale=1]{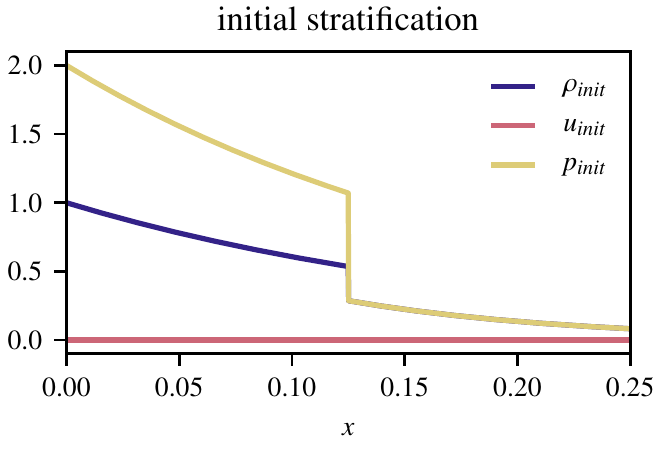}
	\hfill
	\includegraphics[scale=1]{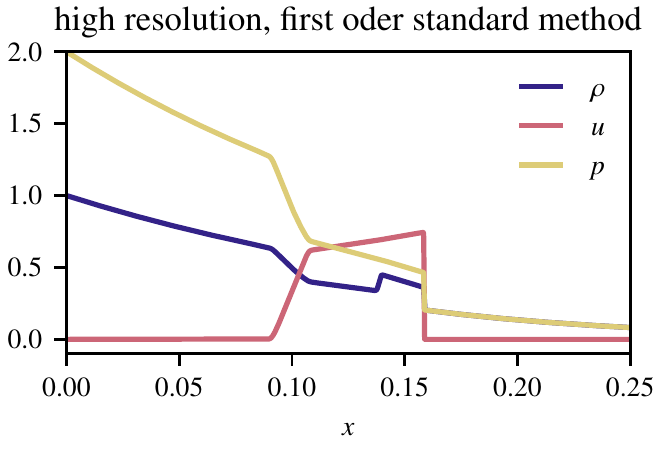}
	\\
	\includegraphics[scale=1]{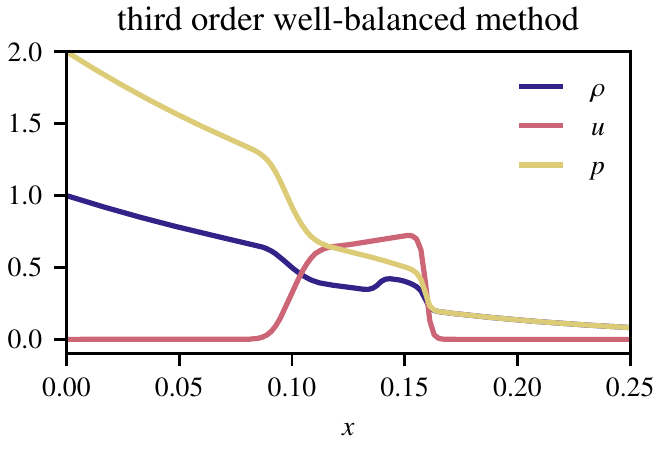}
	\hfill
	\includegraphics[scale=1]{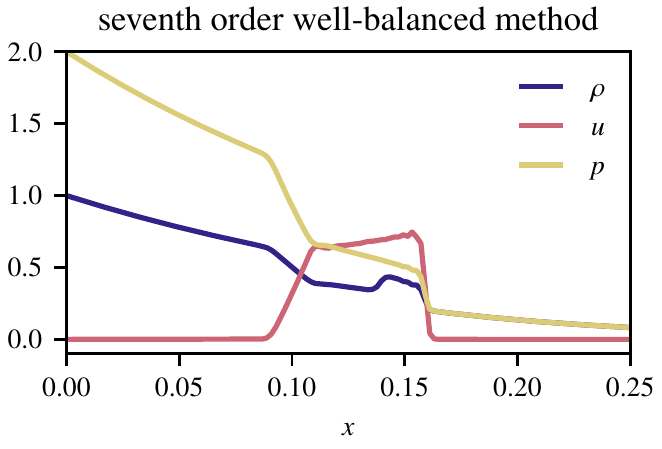}
	\caption{%
			\label{fig:robustness}Riemann problem on an isothermal hydrostatic solution test case from \cref{sec:robustness}. The left top panel shows the initial condition \cref{eq:hystatWithRiemann}. The other panels show numerical results at final time $t=0.02$ as described in the text.
		}
\end{figure}

To test the robustness of our well-balanced methods in combination with CWENO reconstruction we use the initial data
\begin{align}
\label{eq:hystatWithRiemann}
\rho(x)&:=\left\{\begin{matrix}
\exp\left(-\half\phi(x)\right) & \text{if $x<0.125$},\\
\exp(-\phi(x)) & \text{if $x\geq 0.125$},
\end{matrix}\right.
&
 p(x)&:=\left\{\begin{matrix}
2\exp\left(-\half\phi(x)\right) & \text{if $x<0.125$},\\
\exp(-\phi(x)) & \text{if $x\geq 0.125$},
\end{matrix}\right.
&
u(x)&:=0
\end{align}
with $\phi(x):=-10 x$.
\cref{eq:hystatWithRiemann} describes a \revb{piecewise} isothermal hydrostatic solution with a jump, which includes all three waves of the Euler equations. We set these initial data on the domain $[0,0.25]$ and evolve them to the final time $t=0.02$ using our third and seventh order well-balanced method with CWENO reconstruction and Roe's approximate Riemann solver on $128$ grid cells. As \wbref solution for the well-balanced methods we choose 
\begin{align}
\tilde\rho(x)&:=
\exp(-\phi(x))
&
\tilde p(x)&:=
\exp(-\phi(x))
&
\tilde u(x)&:=0.
\end{align}
The results at final time are presented in \cref{fig:robustness}. As a reference we use a result obtained with a first order standard method on 8192 grid cells (top right panel). Neither the third order method (bottom left panel) nor the seventh order method (bottom right panel) show significant oscillations. The wave structure is captured correctly by both methods.

\subsection{2-D numerically approximated hydrostatic solution of the Euler system with gravity}
\label{sec:integrated_hystat}
In stellar astrophysics applications, the hydrostatic state of the star can often be given in a discrete form. In this test, we will show that our well-balanced method can be used if the \wbref solution is given in the form of discrete data in a table.

\begin{figure}
	\centering
	\includegraphics[scale=1.]{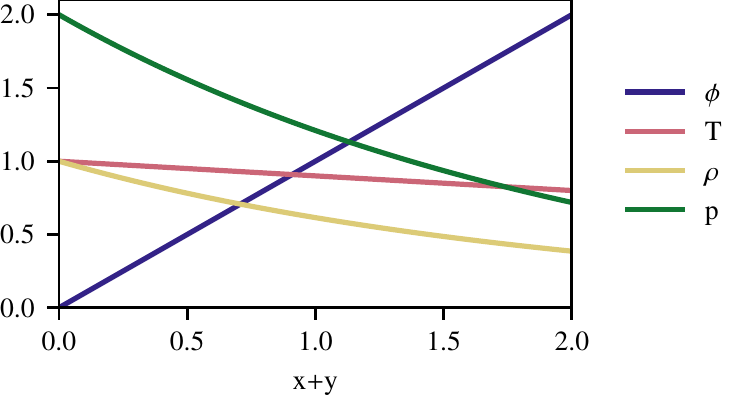}
	\caption{\label{fig:integrated_hystat} Data of the numerically integrated hydrostatic solution used in \cref{sec:integrated_hystat}.}
\end{figure}

\revb{
	The thermodynamical quantities shall be related by the EoS for an ideal gas with radiation pressure, which is given by\cite{Chandrasekhar1958}
	\begin{equation}
	p = \rho T + T^4,
	\label{eq:radiative_eos}
	\end{equation}
	where the temperature $T$ is defined implicitly via
	\begin{equation}
	\epsilon = \frac{T}{\gamma-1} + \frac{3}{\rho} T^4.
	\label{eq:radiative_epsilon}
	\end{equation}
	We assume the following data is given. Let the gravitational potential be $\phi(\x)=\phi(x,y)=x+y$ and the hydrostatic temperature profile is $\bar T(\x)=1-0.1\phi(\x)$.
}
Using Chebfun \cite{Driscoll2014} in the numerical software MATLAB we solve the 1-D hydrostatic equation and EoS for density and pressure corresponding to the given temperature profile. The result is shown in \cref{fig:integrated_hystat}.The data are stored as point values on a fine grid (10,000 data points). The data are set on the 2-D grid using the procedure (a) from \cref{sec:discrete_reference_example1}
We use a $64\times64$ grid to evolve the hydrostatic initial condition to the final time $t=2$.
For the conversion between pressure and internal energy we use Newton's method to solve for the temperature. The $L^1$-errors at final time are shown in \cref{tab:integrated_hystat}. In all tests using the well-balanced modification, there is no error at the final time. 

\begin{table}
	\centering
	\caption{
		\label{tab:integrated_hystat}
		$L^1$-errors in total energy for the numerically integrated hydrostatic solution of the Euler equations with radiation pressure after time $t=2$ computed using different methods. The setup is described in \cref{sec:integrated_hystat}.
	}
		\begin{tabular}{| c | c | c | c | c | c | c |}
			\hline
			$N$ & Std-O1 & WB-O1 & Std-O2 & WB-O2 & Std-O3 & WB-O3 \\
			\hline
			64  & 5.10e-03 & 0.00e+00 & 4.38e-05 & 0.00e+00 & 1.36e-07 & 0.00e+00\\
			\hline
		\end{tabular}
\end{table}
\subsection{Double Gresho vortex}
\label{sec:double_gresho}
\begin{figure}
	\centering
	\includegraphics[scale=1]{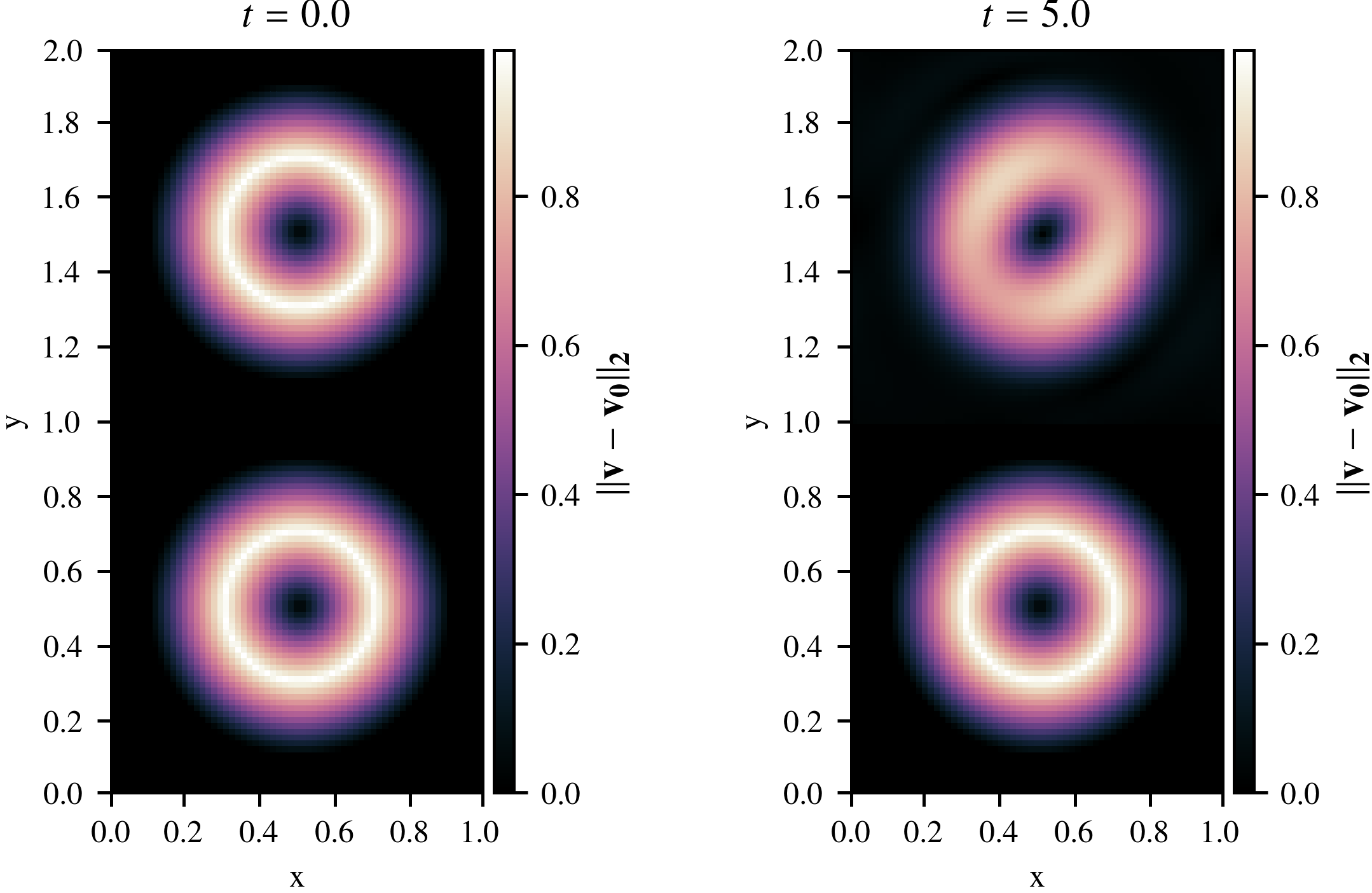}
	\caption{Illustration for the double Gresho vortex test from \cref{sec:double_gresho}. The absolute velocity after subtraction of the constant advection velocity is shown for the initial (left panel) and final (right panel) time. The vortex which is included in the \wbref solution \reva{(which is the bottom vortex in both panels)} is preserved while the other one is diffused and deformed.}
	\label{fig:double_gresho}
\end{figure}

In this test we use a vortex for homogeneous 2-D Euler equations first introduced in \cite{Gresho1990}. The pressure and angular velocity of this vortex in dependence of the distance to the center are given by
\begin{align}
\hat u(r) &= \begin{cases}
5 r, 	& \text{if } 0\leq r<0.2,\\
2-5 r,	& \text{if } 0.2\leq r<0.4\\
0,		& \text{if } 0.4 \leq r,
\end{cases}
\nonumber\\
\hat p(r) &= \begin{cases}
5+\frac{25}2r^2, 	& \text{if } 0\leq r<0.2,\\
9-4\ln(0.2)+\frac{25}{2}r^2-20r+4\ln(r),	& \text{if } 0.2\leq r<0.4\\
3+4\ln(2),		& \text{if } 0.4 \leq r.
\end{cases}
\end{align}
The radial velocity is zero and the density is $\rho\equiv1$.
In our test we set up the domain $[0,1]\times[0,2]$ with two Gresho vortices centered at $(0.5,0.5)$ and $(0.5,1.5)$ respectively. The vortices are advected with the velocity $\vec v_0 = (u_0,v_0)^T=(0.2,0.4)^T$ and the boundaries are periodic. At time $t=5$ the exact solution of this initial data equals the initial setup. We apply our well-balanced method on a $64\times128$ grid to evolve the initial condition up to final time $t=5$. We use Roe's numerical flux functions and a linear reconstruction. Only the vortex initially (and finally) centered at $(0.5,0.5)$ is included in the \wbref solution. The result is illustrated in \cref{fig:double_gresho}.

\subsection{2-D Euler wave in gravitational field}
\label{sec:xingshu}
To demonstrate that we can follow time-dependent solutions exactly with our method we use a problem from \cite{Xing2013} and \cite{Chandrashekar2017} which involves a known exact solution of the 2-D Euler equations with gravity given by
\begin{align}
\tilde\rho(t,x,y) &= 1+\frac{1}{5} \sin(\pi(x+y-t(u_0+v_0))), \qquad \tilde u(t,x,y) = u_0, \qquad \tilde v(t,x,y) = v_0,\nonumber\\
\tilde p(t,x,y) &= p_0 + t(u_0+v_0)-x-y+\frac{1}{5\pi}\cos(\pi(x+y-t(u_0+v_0))).
\label{eq:xingshu}
\end{align}
The gravitational potential is $\phi(\x)=x+y$, the EoS is the ideal gas EoS. In accordance to \cite{Xing2013} and \cite{Chandrashekar2017} we choose $u_0=v_0=1$, $p_0=4.5$ on the domain $[0,1]^2$. \revb{We use the first, second, and third order accurate well-balanced method to evolve the initial data with $t=0$ to a final time $t=0.1$ on a $64\times64$ Cartesian grid and the second order well-balanced method on a $64\times64$ polar grid.
The $L^1$-error in every component of the state vector is exactly zero in each of the tests. We omit showing a table since it does not provide additional insight. 
}
\subsection{Perturbation on the 2-D Euler wave in gravitational field}
\label{sec:xingshu_pert}
\begin{table}
	\centering
	\caption{
		\label{tab:xingshu_pert}
		$L^1$-errors and convergence rates in total energy for different pressure perturbations ($\eta=0.1,10^{-5}$) on the wave in a gravitational field solution of the 2-D Euler equations after time $t=0.1$. The third order standard and well-balanced method are used. The setup is described in \cref{sec:xingshu_pert}.}

	\begin{tabular}{| c | c c | c c | c c | c c |}
		\hline
		\multirow{3}{*}{grid cells} & \multicolumn{4}{c|}{ $\eta=0.1$} & \multicolumn{4}{c|}{$\eta=10^{-5}$}\\
		\cline{2-9}
		& \multicolumn{2}{c|}{Std-O3} & \multicolumn{2}{c|}{WB-O3} & \multicolumn{2}{c|}{Std-O3} & \multicolumn{2}{c|}{WB-O3} \\
		& $E$ error & rate& $E$ error & rate& $E$ error & rate& $E$ error & rate \\
		\hline
		$64\times64$  & 3.96e-04 & \multirow{2}{*}{2.8} & 3.84e-04 & \multirow{2}{*}{2.8} &4.38e-05&\multirow{2}{*}{2.7}& 9.38e-08 & \multirow{2}{*}{2.7}
		\\
		$128\times128$ &  5.81e-05& \multirow{2}{*}{3.0} & 5.65e-05 & \multirow{2}{*}{3.0} &6.68e-06&\multirow{2}{*}{3.0}& 1.49e-08 & \multirow{2}{*}{2.9}
		\\
		$256\times256$ & 7.50e-06 & \multirow{2}{*}{3.1} & 7.30e-06 & \multirow{2}{*}{3.1} &7.72e-07&\multirow{2}{*}{3.1}& 1.93e-09 & \multirow{2}{*}{3.2}
		\\
		$512\times512$ & 8.50e-07 &                      & 8.30e-07 &                      &5.46e-08&& 2.17e-10 &                      
		\\
		\hline
	\end{tabular}

\end{table}

In this test we want to verify the order of accuracy for perturbations to time-dependent \wbref solutions if the well-balanced method is used. For this we use the initial setup from \cref{eq:xingshu} and add a pressure perturbation:
\begin{align}
\rho(t=0,x,y)&=\tilde{\rho}(t=0,x,y),\; 
u(t=0,x,y)    =\tilde u(t=0,x,y), \;  
v(t=0,x,y)    =\tilde v(t=0,x,y), \nonumber\\ 
p(t=0,x,y)   &=\tilde{p}(t=0,x,y)+\eta \exp\left(-100\left(\left(x-\half\right)^2+\left(y-\half\right)^2\right)\right).
\end{align} 
We evolve these initial data to time $t=0.1$ using the third order standard and well-balanced method with $\eta=0.1$. The $L^1$ errors and corresponding convergence rates are presented in \cref{tab:xingshu_pert}. As reference solution for determining the error we use a numerically approximated solution computed using the third order standard method on a $1024^2$ grid. In this test we use exact boundary conditions for the standard method, which means that we evaluate the states in the ghost cells at any time from \cref{eq:xingshu}. We see third order convergence for both methods. However, there seems to be no significant benefit from using the well-balanced method in this test. The choice of $\eta=0.1$ leads to a large discretization error in the perturbation which seems to dominate the total error. \revb{Choosing a large perturbation} is necessary since we use a solution computed from the standard method as a reference to compute the errors. For smaller perturbations the standard method fails to provide a sufficiently accurate reference solution. To yet show the improved accuracy of the well-balanced modification we add a convergence test with a small perturbation of $\eta=10^{-5}$ for which a sufficiently accurate reference solution is produced using the third order well-balanced method on a $1024^2$ grid. The errors and convergence rates for the third order accurate standard and well-balanced methods can also be seen in \cref{tab:xingshu_pert} and it gets evident that the well-balanced method is significantly more accurate on the small perturbation.

\subsection{2-D Keplerian disk}
\label{sec:keplerian}
\begin{figure}[h!]
	\centering
	\includegraphics[scale=1]{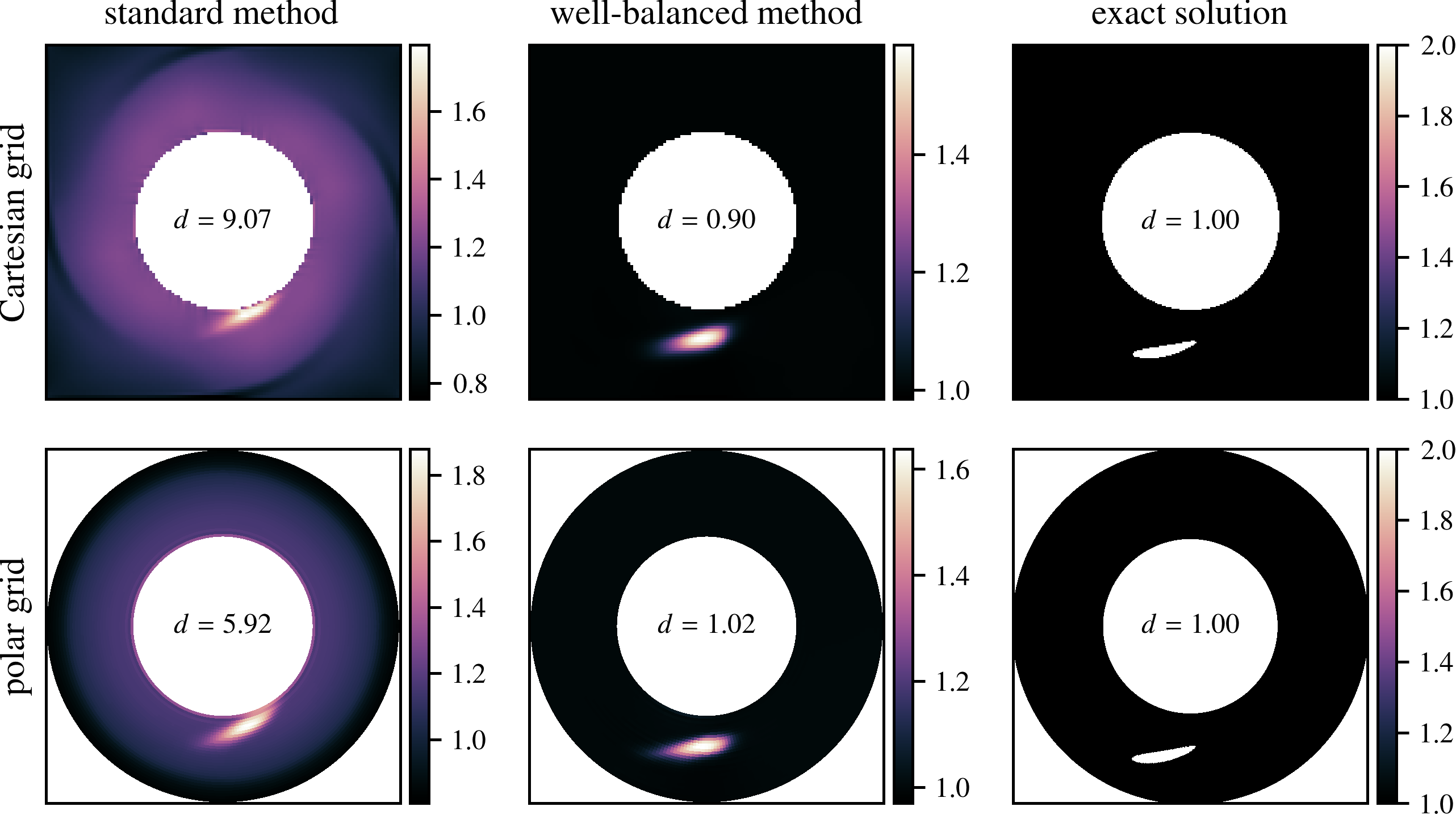}
	\caption{\label{figs:keplerian}Mass advection on a Keplerian disk using second order standard (left column) and well-balanced (center column) methods on a Cartesian $128\times128$ grid (top row) and a polar $32\times256$ grid (bottom row). The setup is given in \cref{sec:keplerian}. Results can be compared to the exact solutions shown in the right column. The meaning of $d$ is described in the text. The domain shown in each panel is $[-2,2]^2$. The $x$-coordinate increases to the right, the $y$-coordinate increases to the top.}
\end{figure}
Consider a stationary solution given by \cite{Gaburro2018b}
\begin{align}
\tilde{\rho}&\equiv 1,& \tilde u(x,y)&=-\sin(\alpha(x,y))\sqrt{\frac{G m_S}{r(x,y)}},& \tilde v(x,y)&=\cos(\alpha(x,y))\sqrt{\frac{G m_S}{r(x,y)}},& \tilde p&\equiv 1
\label{eq:keplerian}
\end{align}
with the gravitational potential $\phi(r)=-\frac{G m_s}r$ and $r=\sqrt{x^2+y^2}$, $\alpha=\arctan(\frac y x)$, $G=m_s=1$. We use the initial conditions
\begin{equation}
\rho(x,y) = \begin{cases}2&\text{if }(x-1.2)^2+(y-1.0)^2<0.15^2\\\tilde{\rho}&\text{else}\end{cases}
\end{equation}
and $u=\tilde u, v=\tilde v,p=\tilde p$ on the domain $[-2,2]\times[-2,2]$. We chose the domain such that we omit the singularity in the velocity at $(x,y)=(0,0)$. In \cref{figs:keplerian} results of numerical tests are illustrated. The second order standard and well-balanced methods are applied on a polar grid with $32\times256$ cells and a Cartesian grid with $128\times128$ cells. In the Cartesian grid we take out the center with $r<1$ using Dirichlet boundary conditions. We also use Dirichlet boundary conditions at all outer boundaries. Since there is a discontinuity in the initial setup we apply a minmod slope limiter \revb{to the linear reconstruction of the conserved variables in the standard method or the deviations in conserved variables in the well-balanced method}. The density at time $t=2.5$ for each simulation is shown in \cref{figs:keplerian} together with the exact solution. Since this is an purely advective problem and there is no radial component to the velocity, the quantity $\|(\rho-1)r\|_1$, which describes the average distance of the density perturbation to the center, is conserved for all time in the exact solution. For our simulations we measure the quantity $d = \|(\rho(t=2.5)-1)r\|_1/\|(\rho(t=0)-1)r\|_1$ as a measure of the quality of the numerical solutions. For the exact solution we have $d=1$ for all time. The values of $d$ are shown in the center of the plots in \cref{figs:keplerian}.

In the tests with the standard method we see discretization errors in the Keplerian disk solution \cref{eq:keplerian}. This introduces radial velocities, the advection of the spot of increased density has a component towards the center.
In the tests using our well-balanced methods the result is free of discretization errors in the Keplerian disk solution \cref{eq:keplerian}. The advection is more accurate, the only errors are diffusion errors. The polar grid is more suitable for this test problem, since it is adapted to the radial geometry. The test using our well-balanced method on the Cartesian grid is more diffusive than the one on the polar grid, yet we see that the well-balanced modification improves the result significantly on both grids.

\subsection{Stationary MHD vortex - long time}
\label{sec:mhd_vortex_long_time}
We consider the following exact solution of the homogeneous 2-D ideal MHD equations:
\begin{align}
\hat x &= x - t u_0,&
\hat y &= y - t v_0,&
r^2 &=\hat x^2+\hat y^2, \nonumber\\
u&=u_0-k_pe^{\frac{1-r^2}{2}} \hat y,&
v&=v_0+k_pe^{\frac{1-r^2}{2}} \hat x,&
\rho&=1,\nonumber\\
B_x&=-m_p e^{\frac{1-r^2}{2}} \hat y,&
B_y&= m_p e^{\frac{1-r^2}{2}} \hat x,&
p&= 1+ \left(\frac{m_p^2}{2} (1-r^2) -\frac{k_p^2}{2}\right)e^{1-r^2}.
\end{align}
This setup describes a stationary vortex which is advected through the domain with the velocity $(u_0,v_0)$. The domain is $[-5,5]\times[-5,5]$.
One vortex turnover-time is $t_\text{turnover}=\frac{2\pi}{\sqrt{e}k_p}\approx\frac{3.81}{k_p}$. In a first test we set $m_p=k_p=0.1$, $u_0=v_0=0$ and run the test up to $t=100 t_\text{turnover}$ on a $32\times32$ grid. We use the well-balanced method and the \wbref solution equals the initial data. The numerical error at final time compared to the initial setup is exactly zero in all conservative variables.
\subsection{Stationary MHD vortex - order of accuracy}
\label{sec:mhd_vortex}

In a second test with the vortex from \cref{sec:mhd_vortex_long_time} we want to see if the well-balanced method converges as expected, even if the \wbref solution which is set deviates from the actual solution over time. For that we set $m_p=k_p=0.1$, $u_0=v_0=0$ in the initial condition. As \wbref solution we use the same vortex but with $u_0=v_0=1$.
In \cref{tab:mhd_vortex} the $L^1$ errors and rates in energy at final time $t=0.2$ are presented for the formally first, second, and third order accurate well-balanced method. We see that even if the \wbref solution moves away from the actual solution over time the method is consistent with the expected order of accuracy.
\begin{table}
	\centering
	\caption{
		\label{tab:mhd_vortex}
		$L^1$-errors and convergence rates in total energy for a 2-D MHD stationary vortex after time $t=0.2$. Different methods are used. In the well-balanced method, a \wbref solution is chosen, which deviates from the actual solution over time. The setup is described in~\cref{sec:mhd_vortex}.
	}
	\begin{tabular}{| c | c c | c c | c c | c c|}
			\hline
			\multirow{2}{*}{grid cells} & \multicolumn{2}{c|}{WB-O1} & \multicolumn{2}{c|}{WB-O2} & \multicolumn{2}{c|}{WB-O3} & \multicolumn{2}{c|}{Std-03}   \\
			& $E$ error & rate& $E$ error & rate& $E$ error & rate& $E$ error & rate \\
			\hline
			$128 \times 128$  & 2.42e-02 & \multirow{2}{*}{1.0} & 1.87e-03 & \multirow{2}{*}{1.9} & 8.46e-05 & \multirow{2}{*}{3.0} & 8.46e-06 & \multirow{2}{*}{3.3}
			\\
			$256 \times 256$  & 1.23e-02 & \multirow{2}{*}{1.0} & 5.11e-04 & \multirow{2}{*}{1.9} & 1.08e-05 & \multirow{2}{*}{3.0} & 8.38e-07 & \multirow{2}{*}{3.1}
			\\
			$512 \times 512$  & 6.17e-03 &                      & 1.34e-04 &                      & 1.36e-06 &                      & 9.67e-08 &                      
			\\
			\hline
	\end{tabular}
\end{table}

	\subsection{Stationary MHD vortex - numerical \wbref solution}
	\label{sec:mhd_vortex_numref}
\begin{figure}
		\includegraphics[width=0.98\textwidth]{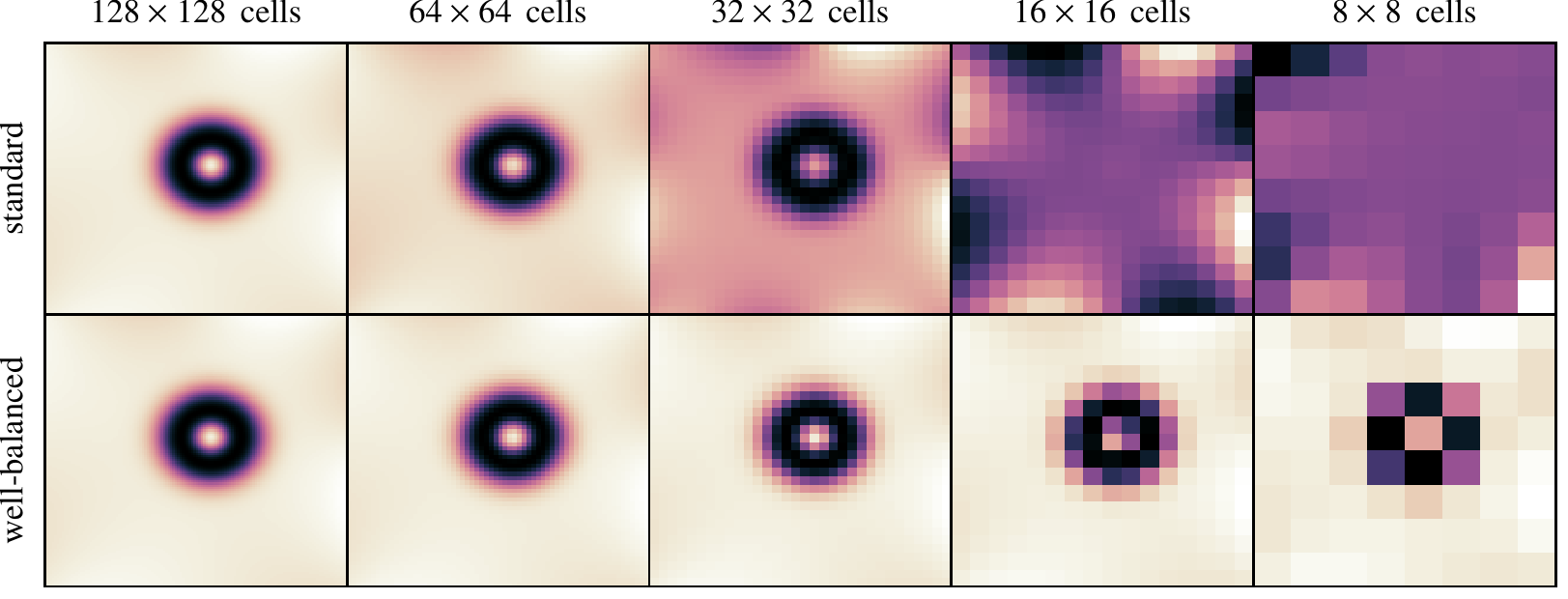}
		\caption{Pressure of the moving stationary MHD vortex as described in \cref{sec:mhd_vortex_numref} at the final time. The \wbref solution for the well-balanced method is the numerical solution computed with the standard method on a $128\times 128$ cells grid (upper left panel). In the upper panels the standard method is used, in the lower panels the well-balanced method is used. Different columns correspond to different resolutions.\label{fig:mhd_vortex_numref}}
	\end{figure}

	In this test we present a simple application of the method described in \cref{sec:discrete_reference_example2}. Again, we use the stationary MHD vortex test case described in \cref{sec:mhd_vortex_long_time}. The parameters are $k_p=m_p=0.1$ and $u_0=v_0=0.1$, the final time is $t_\text{final}=5$.
	First, we compute a \wbref solution using our third order non-well-balanced method with $128\times128$ grid cells with parabolic extrapolation boundary conditions. Every time-step is stored. This is then used in our well-balanced method as described in \cref{sec:discrete_reference_example2} using third order interpolation in time and in space. The resulting pressure for well-balanced and standard methods on different Cartesian meshes is shown in \cref{fig:mhd_vortex_numref}. All methods use CWENO3 reconstruction and parabolic extrapolation boundary conditions.
	On the $128\times128$ grid the solutions for the well-balanced and non-well-balanced method are exactly the same, since the solution from the standard method is used as \wbref solution in the well-balanced method.
	For smaller resolutions the standard method is too diffusive to resolve the vortex. The quality of the results obtained with the well-balanced method is the same for all resolutions, since all of them use the same $128\times128$ simulation as \wbref solution. 

\section{Computational cost of the modification}
\label{sec:efficiency}

Well-balanced methods are constructed to improve the accuracy with which solutions of balance laws are approximated. In the previous section we have shown that usage of our well-balanced modification can improve the accuracy of a simulation significantly. On the other hand, an increase of computational effort can countervail the gain in accuracy, if it is too high. In this section we will compare the computation times of simulations using our well-balanced modification to simulations using the corresponding standard method and show that the increase in CPU time is moderate.
\subsection{The procedure}
\label{sec:efficiency_procedure}

To compare the methods, we will run tests with different setups and grid resolutions using a standard method and the corresponding well-balanced method. We use the simple python code described in \cref{sec:numerical_tests} on a single CPU. Each test is repeated 20 times and the wall clock times are measured. We compute the average wall-clock time and standard deviations of the single runs for every test. The ratio of average wall clock time for the well-balanced compared to the standard method is visualized depending on the grid resolution.

\revb{
Note that we use final times that are significantly smaller than the final times used in the corresponding tests above. The reason for this is that after some time the solutions obtain with and without well-balancing differ. This results in different sizes for the time steps which significantly influences the wall-clock time necessary to reach the final time. Since we aim to compare the efficiency of the methods without taking the quality of the solution into account, we use reduced final times.}
\subsection{The tests}
\begin{figure}
	\centering
	\includegraphics[scale=1]{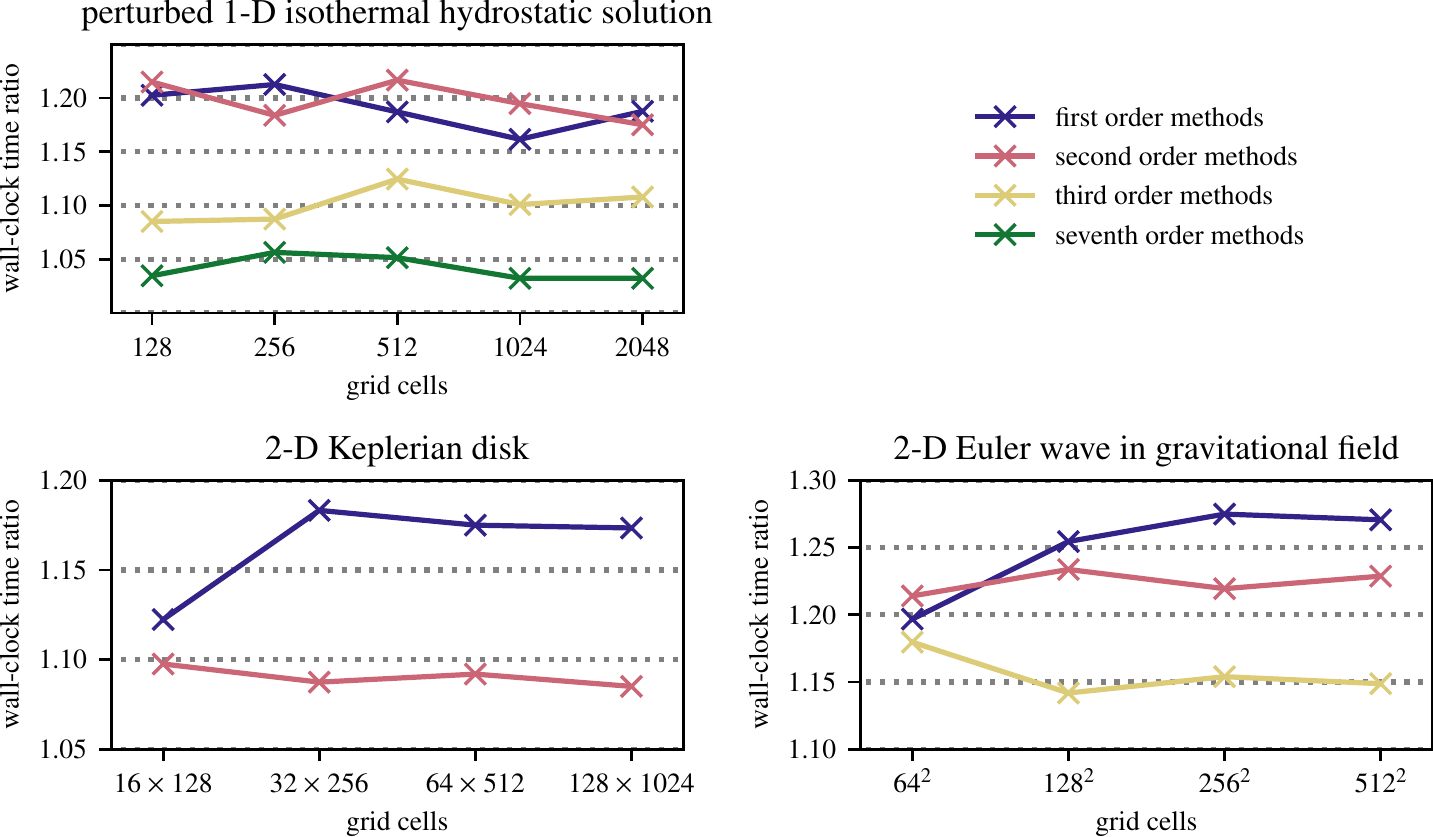}
	\caption{Ratio of the average wall-clock times for the well-balanced and the standard method. Values and errors are determined as described in \cref{sec:efficiency_procedure}.\emph{Top left}: One-dimensional test case with static \wbref solution from \cref{sec:1d_perturbation} \emph{Bottom left}: Two-dimensional test case with stationary \wbref solution from \cref{sec:keplerian} on a polar grid. \emph{Bottom right}: Two-dimensional test case with time-dependent \wbref solution from \cref{sec:xingshu_pert} on a Cartesian grid.}
	\label{fig:efficiency}
\end{figure}

To compare the runtimes we use test setups described in \cref{sec:numerical_tests} using exactly the same methods. In a first test, we use the perturbed one-dimensional isothermal solution described in \cref{sec:1d_perturbation} with a final time $t=0.02$. The first, second, third, and seventh order 1-D methods are applied. The ratio of wall clock times for the tests with and without the well-balanced modification can be seen in the top right panel of \cref{fig:efficiency}.  The \wbref solution $\tilde q$ used in the well-balanced modification is constant in time in this test case. It is hence computed once and stored in an array. From the figure we see that we can expect an increase in CPU time of about $20\%$ for the first and second order accurate methods. The difference in runtime reduces significantly for methods with higher order of accuracy. For the seventh order methods the increase in wall-clock time is only about $5\%$.

As a second test setup we choose the Keplerian disk from \cref{sec:keplerian}. We evolve it to time $t=0.01$ with a first and second order method on a polar grid. As in the previous test, the \wbref solution is time-independent and thus only computed once. The wall-clock time ratios are visualized in the bottom left panel of \cref{fig:efficiency}. We observe an increase in CPU time of less than $20\%$ when using the first order well-balanced method and less than $10\%$ when using the second order well-balanced method.

To also test the increase of CPU time consumption for a simulation in which the \wbref solution is time dependent, we use the Euler wave in a gravitational field with perturbation from \cref{sec:xingshu_pert}. We evolve the solution up to the final time $t=0.01$ with each method. In this test, the \wbref solution is computed from a function every time it is used (which happens in every intermediate step). The result of these tests can be seen in the bottom right panel of \cref{fig:efficiency}. We see an increase in CPU time of less than $30\%$ if the well-balanced method is used. Note that, again, the wall-clock time ratio is smaller for methods with higher order accuracy. For the third order well-balanced method we only observe an increase of about $15\%$ in wall-clock time.

\section{Summary and conclusions}

We introduced a new general framework for the construction of well-balanced finite volume methods for hyperbolic balance laws. A standard finite volume method is modified such that it evolves the deviation from a \wbref solution instead of the actual solution. This makes the scheme exact on the \wbref solution. The finite volume method can include any consistent reconstruction, numerical flux function, interface quadrature, source term discretization, and ODE solver for time discretization. Thus, it can be arbitrarily high order accurate and the method can be defined on any computational grid geometry. One can view our method as a high order extension of \cite{Berberich2019} and \cite{Klingenberg2019} to all known solutions of all hyperbolic balance laws.

In numerical tests with Euler and MHD equations on different grids we could verify that the method can successfully be applied to exactly maintain static and stationary solutions or even follow time-dependent solutions. For that, the solution has to be known either analytically or in the form of discrete data. The latter case is especially interesting for complex applications like stellar astrophysics, where static states of the Euler equations with gravity can be obtained numerically but only in few cases analytically. Also, for the case of differentially rotating stars, \revb{\ie stars that are in a rotating stationary state with an angular velocity depending on the distance to the axis of rotation}, our method can be applied for well-balancing since it can include non-zero velocities in stationary states. High order accuracy has been confirmed in numerical experiments. Also, in a series of numerical tests we have shown that the increase in computational time is moderate. 

\reva{
The proposed well-balanced method can be easily implemented in existing finite volume codes with minimal effort. However, there are applications in which the well-balanced solution is not known beforehand and the current method cannot be applied. In that case, another  well-balanced methods existing in the literature for the balance law under consideration has to be applied, in case such a scheme exists. In all other cases, in which the well-balanced solution is known, our simple \revb{framework} can be applied to obtain the well-balanced property.
}
\section*{Acknowledgments}
The research of Jonas Berberich is supported by the Klaus Tschira Foundation. Praveen Chandrashekar would like to acknowledge the support received from Airbus Foundation Chair on Mathematics of Complex Systems at TIFR-CAM, Bangalore and Department of Atomic Energy, Government of India, under  project no.  12-R\&D-TFR-5.01-0520.

\bibliographystyle{model1-num-names}
\bibliography{library}
\appendix
\section{Details of the applied finite volume schemes}
\label{appendix:detail}
\renewcommand{\P}{\ensuremath{\mathcal P}\xspace}
\newcommand{\U}{\vec U}
We use structured grids in all the numerical tests and hence in the description of the details, we restrict ourselves to structured grids. Some parts of the scheme, such as the reconstruction methods, are applied to $\Q$ in the standard method but to $\Delta\Q$ in the well-balanced method. We will denote the states with $\U$; depending on the method we have $\U=\Q$ or $\U=\Delta\Q$. Analogously, we will use $\vec u$ to denote $\q$ or $\Delta \q$.
\subsection{Curvilinear coordinates}
\label{appendix:curvilinear}
\newcommand{\ud}{\textrm{d}}
\newcommand{\dd}[2]{\frac{\ud #1}{\ud #2}}
\newcommand{\df}[2]{\frac{\partial #1}{\partial #2}}
\newcommand{\con}{\bm{q}}
\newcommand{\ev}{\bm{q}}
\newcommand{\ee}{\textrm{e}}

\newcommand{\nfl}{\hat{\fl}}
\newcommand{\bs}{\bm{s}}
\newcommand{\avg}[1]{\bar{#1}}
\newcommand{\lavg}[1]{\hat{#1}}
\newcommand{\mm}[1]{\textrm{mm}\left( #1 \right)}
\newcommand{\hT}{\hat{T}}
\newcommand{\hbeta}{\hat{\beta}}
\newcommand{\hcon}{\bm{Q}}
\newcommand{\hfl}{\bm{F}}
\newcommand{\hgl}{\bm{G}}
\newcommand{\hbs}{\bm{S}}

\newcommand{\nhfl}{\hat{\hfl}}
\newcommand{\nhgl}{\hat{\hgl}}

\newcommand{\clr}[2]{\color{#1} #2 \color{black}}

\newcommand{\fl}{\bm{f}}
\newcommand{\gl}{\bm{g}}
\newcommand{\bw}{\bm{w}}
We define a 2-D curvilinear coordinate system. The coordinates in physical space are $\x=(x,y)$, the coordinates in computational space are $\boldsymbol\xi=(\xi,\eta)$.
The $(i,j)$-th cell is denoted $\Omega_{i,j}$ in the physical space and by $\hat{\Omega}_{i,j}$ in the computational space. 
We can rewrite \cref{eq:eul2d_cartesian} in the computational coordinates as
\begin{equation}
J \dt\q + \partial_\xi\hat\f + \partial_\eta\hat{\g}=J \s,
\label{eq:euler_curvilinear}
\end{equation}
where
\begin{align}
J:&=\textrm{det}\begin{pmatrix}
	\DD x\xi & \DD x\eta \\
	\DD y\xi & \DD y\eta
\end{pmatrix},&
\hat{\f} :&= J\left(\DD \xi x\f + \DD\xi y\g\right), &
\hat{\g} :&= J\left(\DD \eta x\f + \DD\eta y\g\right).
\end{align}
To solve \cref{eq:eul2d_cartesian} on the curvilinear physical grid, we can now solve \cref{eq:euler_curvilinear} on a Cartesian grid. We construct the grid from the nodes and approximate the derivatives of the coordinate \reva{transformation} using central differences on the nodal coordinates. This implementation of curvilinear grids restricts the scheme to second order accuracy. We can achieve higher order accuracy only on Cartesian grids. More details on the finite volume method on a curvilinear mesh can be found in \cite{Berberich2019}.
\\\noindent\emph{Polar grid}: The polar grid can be defined by the function
\begin{equation}
	\x(\boldsymbol\xi) := 
	\begin{pmatrix}
		\xi\sin(\eta)\\\xi\cos(\eta)
	\end{pmatrix}
\end{equation}
for $\xi>0$, $\eta\in\left[0,2\pi\right)$. Note, that this functions can not be inverted at $\xi=0$, i.e.\ $\x=0$. Hence, the origin in physical coordinates has to be omitted, when this grid is used.

\subsection{Implementation of boundary conditions}
\label{appendix:boundaries}

Boundary conditions are implemented on the structured grid by using ghost cells to artificially extend the computational domain. In this section we assume the physical domain is in the cells $\Omega_{ij}$, $(i,j)\in\{0,\dots,N\}\times\{0,\dots,M\}$. The necessary amount of ghost cells depends on the stencil of the method, especially on the quadrature and reconstruction (plus one row of ghost cells for the flux computation at the boundary).
The values in the ghost cells are set after each intermediate Runge--Kutta step. In the following, the used boundary conditions are presented for 2-D grids. If there is no description for the 1-D method, the method is reduced to 1-D in the trivial way.

%
\noindent\emph{Constant extrapolation}:
The constant extrapolation boundary conditions are suitable to be used with a first oder accurate method. They are obtained by setting
\begin{equation}
	\U_{ij}=\U_{kl}\qquad\text{with}\qquad k=\min(N,\max(0,i)),\quad l=\min(M,\max(0,j))
\end{equation}
in the ghost cells.

\noindent\emph{Linear extrapolation}:
In combination with a second order accurate method on a Cartesian grid, the following linear extrapolation to the ghost cells can be used as boundary conditions: We set
\begin{align}
	\U_{-k,j} &= (1+k)\U_{0,j}-k\U_{1,j},&
	\U_{i,-k} &= (1+k)\U_{0,M}-k\U_{i,1},\nonumber\\
	\U_{N+k,j} &= (1+k)\U_{N,j}-k\U_{N-1,j},&
	\U_{i,M+k} &= (1+k)\U_{i,M}-k\U_{i,M-1}
\end{align}
for $(i,j)\in\{0,\dots,N\}\times\{0,\dots,M\}$, $k=1,2$. The diagonal ghost cells ($\Omega_{N+1,M+1}$, etc.) are not needed for the second order scheme used in our tests.

\noindent\emph{Parabolic extrapolation}:
Parabolic extrapolation to the ghost cells is suitable to be used in combination with a third order accurate method. On a 1-D equidistant grid we set the values in the ghost cells to
\begin{align}
	\U_{N+1}&=3 \U_N-3 \U_{N-1}+\U_{N-2},&
	\U_{N+2}&=6 \U_N-8 \U_{N-1}+2\U_{N-2},\nonumber\\
	\U_{-1}&=3 \U_0-3 \U_{1}+\U_{2},&
	\U_{-2}&=6 \U_0-8 \U_{1}+2\U_{2}.
\end{align}
On a two-dimensional grid we have to use a genuine 2-Dimensional parabolic extrapolation to obtain third order accuracy. The basis for that is a 2-D parabolic reconstruction parabola. We use the $P_\text{opt}$ parabola from \cite{Levy2000} which is given by
\begin{multline}
	\vec u(x,y)\approx \frac{\Delta x\Delta y}{24} \Bigg(\frac{12 (x-x_i)^2 (-2 \U_{i,j}+\U_{i+1,j}+\U_{i-1,j})}{\Delta x^2}\\
	+\frac{6 (x-x_i) (y-y_j) (\U_{i+1,j+1}-\U_{i-1,j+1}-\U_{i+1,j-1}+\U_{i-1,j-1})}{\Delta x \Delta y}
	+\frac{12 (\U_{i+1,j}-\U_{i-1,j}) (x-x_i)}{\Delta x}\\
	+\frac{12 (y-y_j)^2 (-2 \U_{i,j}+\U_{i,j+1}+\U_{i,j-1})}{\Delta y^2}+\frac{12 (\U_{i,j+1}-\U_{i,j-1}) (y-y_j)}{\Delta y}\\
	+28 \U_{i,j}-\U_{i+1,j}-\U_{i,j+1}-\U_{i,j-1}-\U_{i-1,j}\Bigg)
\end{multline}
for $(x,y)\in\Omega_{ij}$.
The values in the ghost cells are then computed by integrating the closest possible reconstruction parabola (computed from only values inside the domain) over the ghost cell. This yields
\begin{equation}
	\U_{N+k,j} = \frac{1}{2} \left(k^2 (-2 \U_{N-1,j}+\U_{N,j}+\U_{N-2,j})+k (-4 \U_{N-1,j}+3 \U_{N,j}+\U_{N-2,j})+2 \U_{N,j}\right)
\end{equation}
and correspondingly for $\Omega_{-k,j}$, $\Omega_{i,M+k}$, and  $\Omega_{i,-k}$ with $k=1,2,..$ and $(i,j)\in\{1,N-1\}\times\{1,M-1\}$.
At the upper right edge we obtain
\begin{align}
	\U_{N+k,M+l} = \phantom{+}&\frac{1}{4} k^2 \left(-4 \U_{N-1,M-1}+2 \U_{N,M-1}+2 \U_{N-2,M-1}\right)\nonumber\\
	+&\frac{1}{4} k \left((l+1) (\U_{N,M}-\U_{N-2,M}-\U_{N,M-2}+\U_{N-2,M-2})-8 \U_{N-1,M-1}+6 \U_{N,M-1}+2 \U_{N-2,M-1}\right)\nonumber\\
	+&\frac{1}{4} \Big(-4 (l+1)^2 \U_{N-1,M-1}+2 (l+2) (l+1) \U_{N-1,M}+(l+1) \U_{N,M}\nonumber
	\\&\phantom{\frac{1}{4} \Big(}-(l+1) \U_{N-2,M}
	+2 l (l+1) \U_{N-1,M-2}-(l+1) \U_{N,M-2}+(l+1) \U_{N-2,M-2}+4 \U_{N,M-1}\Big)
\end{align}
for a ghost cell $\Omega_{N+k,M+l}$($k,l\in\{0,1,2,..\}^2\setminus\{0\}^2$) and correspondingly for the other edges.

\subsection{Interface quadrature}

For the two-dimensional third order method it is necessary to apply a quadrature rule to compute the interface flux. For that we use the Gau\ss--Legendre quadrature rule. For that, in \cref{eq:semi-discrete_3d}, we use $M=3$ and the weights
\begin{equation}
	\omega_1=\omega_3=\frac5{18},\qquad \omega_2=\frac49.
\end{equation}
The quadrature points are
\begin{align}
	\x_{i+\half,k,1} &= \sqrt{\frac35} \x_{i+\half,k-\half} + \left(1-\sqrt{\frac35}\right) \x_{i+\half,k},&
	\x_{i+\half,k,2} &= \x_{i+\half,k},&
	\x_{i+\half,k,3} &= \sqrt{\frac35} \x_{i+\half,k+\half} + \left(1-\sqrt{\frac35}\right) \x_{i+\half,k}
	\label{eq:quadrature_points1}
\end{align}
at the interface between the $\Omega_{ik}$ and the $\Omega_{i+i,k}$ control volumes and
\begin{align}
\x_{i,k+\half,1} &= \sqrt{\frac35} \x_{i-\half,k+\half} + \left(1-\sqrt{\frac35}\right) \x_{i,k+\half},&
\x_{i,k+\half,2} &= \x_{i,k+\half},&
\x_{i,k+\half,3} &= \sqrt{\frac35} \x_{i+\half,k+\half} + \left(1-\sqrt{\frac35}\right) \x_{i,k+\half}
\label{eq:quadrature_points2}
\end{align}
at the interface between the $\Omega_{ik}$ and the $\Omega_{i+i,k}$ control volumes. Note the change in notation from \cref{eq:semi-discrete_3d} to \cref{eq:quadrature_points1,eq:quadrature_points2}. In \cref{eq:semi-discrete_3d} the indices $i$ and $k$ are indices of the $i$-th and $k$-tho control volumes. In \cref{eq:quadrature_points1,eq:quadrature_points2} we consider the case of a structured grid and the indices denote the position of the $\Omega_{ik}$ control volume in the grid. The interface is then denoted using half values.

\subsection{Source term discretization}

In some tests we use a gravity source term for Euler equations. The source term component in the momentum equation has to be approximated to sufficiently high order.
\\\noindent\emph{Second order source term discretization}: For the first and second order method, we use the second order accurate source term discretization
\begin{equation}
	-\frac1{\Delta_x}\int_{\Omega_i} \rho(x) g(x) dx \approx - \bar{\rho}\left(\Omega_{i}\right)g(x_{i})
\end{equation}
in the one-dimensional case and
\begin{equation}
-\frac1{|\Omega_{ij}|}\int_{\Omega_{ij}} \rho(\x) \vec g(\x) d\x \approx - \bar{\rho}\left(\Omega_{ij}\right)\vec g(\x_{ij})
\end{equation}
in the two-dimensional case in the momentum equation. The cell-averaged density is denoted $\bar{\rho}$ and the gravitational acceleration $\g=\nabla\phi$ is given exactly at the cell-center.

	\noindent\emph{Third and seventh order source term discretization in 1-D}:
	For the one-dimensional methods with CWENO reconstruction we define 
	$$s^{\rho u}_h(x):=-\rho^\text{rec}(x)g^\text{int}(x),$$
	where $\rho^\text{rec}$ is the density polynomial obtained from the CWENO reconstruction and $g^\text{int}$ is the gravitational acceleration interpolated from the cell centered values to third or seventh order respectively. Since $s^{\rho u}_h$ is a polynomial, the cell-average
	$$S^{\rho u}_i:=\frac1{|\Omega_i|}\int_{\Omega_i}s^{\rho u}_h(x)\,dx$$
	can be computed exactly and is used as source term approximation in the momentum equation in the scheme.
	
	\noindent\emph{Third order source term discretization in 2-D}:
	Similar to the 1-D case we define the vector-valued function
	$$
		(s^{\rho u}_h,s^{\rho v}_h)^T(\x) := -\rho^\text{rec}(\x)\vec g^\text{int}(\x),
	$$
	where $\rho^\text{rec}$ is obtained from the two-dimensional CWENO3 reconstruction and $\g^\text{int}$ is given by
	the interpolation
	\begin{multline*}
		\g^\text{int}(\x) := \g_{ij} + \frac12\left(\g_{i+1,j}-\g_{i-1,j}\right) x + \frac12\left(\g_{i,j+1}-\g_{i,j-1}\right) y 
		+ \frac12\left(\g_{i+1,j}-2\g_{ij}+\g_{i-1,j}\right) x^2\\ +\frac14\left(\g_{i+1,j+1}-\g_{i-1,j+1}-\g_{i+1,j-1}+\g_{i-1,j-1}\right) xy
		+ \frac12\left(\g_{i,j+1}-2\g_{ij}+\g_{i,j-1}\right) y^2.
	\end{multline*}
	for $(x,y)\in\Omega_{ij}$. This polynomial is constructed such that it satisfies $\g^\text{int}(\x_{ij})=\g_{ij}$, $\g^\text{int}(\x_{i\pm1,j\pm1})=\g_{i\pm1,j\pm1},$ and $\g^\text{int}(\x_{i\pm1,j\mp1})=\g_{i\pm1,j\mp1}$ for the cell-centered point values in coordinate direction. The diagonal values $\g_{i+1,j+1}$, $\g_{i-1,j+1}$, $\g_{i+1,j-1}$, and $\g_{i-1,j-1}$ are approximated in the least square sense.
	
	Since $s^{\rho u}_h$ and $s^{\rho v}_h$  are polynomials, the cell-averages
	$$
		S^{\rho u}_{ij}:=\frac1{|\Omega_{ij}|}\int_{\Omega_{ij}}s^{\rho u}_h(x)\,d\x
		\quad\text{and}\quad
		S^{\rho v}_{ij}:=\frac1{|\Omega_{ij}|}\int_{\Omega_{ij}}s^{\rho v}_h(x)\,d\x
	$$
	can be computed exactly and are used as source term approximations in the momentum equations in the scheme.
	
\end{document}